\documentclass[11pt,amssymb]{article}

\usepackage{amsmath,amsthm}

\usepackage{amssymb,latexsym}

\usepackage{enumerate}

\newcommand{\LL}{{\cal L}}
\newcommand{\WW}{{\cal W}}
\newcommand{\XX}{{\cal X}}

\newcommand{\BR}{{\mathbb R}}
\newcommand{\BX}{{\mathbb X}}

\newcommand{\dyw}{\mbox{\rm div}}

\newtheorem{theorem}{\bf Theorem}[subsection]
\newtheorem{proposition}[theorem]{\bf Proposition}
\newtheorem{definition}[theorem]{\bf Definition}
\newtheorem{lemma}[theorem]{\bf Lemma}
\newtheorem{corollary}[theorem]{\bf Corollary}
\newtheorem{remark}[theorem]{\bf Remark}
\newtheorem{example}[theorem]{\bf Example}

\newcommand{\nsubsection}{\setcounter{equation}{0}\subsection}

\begin{document}

\title {On time-dependent functionals of diffusions corresponding
to divergence form operators\footnote{Research supported by the
Polish Minister of Science and Higher Education under Grant N N201
372 436.}}
\author {Tomasz Klimsiak \smallskip\\
{\small Faculty of Mathematics and Computer Science,
Nicolaus Copernicus University} \\
{\small  Chopina 12/18, 87--100 Toru\'n, Poland}\\
{\small e-mail: tomas@mat.uni.torun.pl}}
\date{}
\maketitle \noindent{\small{\bf Abstract:}  We consider processes
of the form $[s,T]\ni t\mapsto u(t,X_t)$, where $(X, P_{s,x})$ is
a multidimensional diffusion corresponding to uniformly elliptic
divergence form operator. We show that if
$u\in{\mathbb{L}}_{2}(0,T;H_{\varrho}^{1})$ with $\frac{\partial
u}{\partial t} \in{\mathbb{L}}_{2}(0,T;H_{\varrho}^{-1})$ then
there is a quasi-continuous version $\tilde u$ of $u$ such that
$\tilde u(t,X_t)$ is a $P_{s,x}$-Dirichlet process for quasi-every
$(s,x)\in[0,T)\times{\mathbb{R}}^d$ with respect to parabolic
capacity and we describe the martingale and the zero-quadratic
variation parts of its decomposition. We give also conditions on
$u$ ensuring that $\tilde u(t,X_t)$ is a semimartingale.}

\nsubsection{Introduction}

In the present paper we study   structure of additive
functionals (AFs for short) of the form $X^u=\{X^u_{s,t}\equiv
u(t,X_t)-u(s,X_s);0\le s\le t\le T\}$, where
$u:Q_T\equiv[0,T]\times{\mathbb{R}}^{d}\rightarrow{\mathbb{R}}$
and $\BX=\{(X,P_{s,x});(s,x)\in Q_T\}$ is a Markov family
corresponding to the operator
\begin{equation}\label{eq1.1}
L_{t}=\frac{1}{2}\sum_{i,j=1}^d\frac{\partial}{\partial
x_{i}}(a_{ij}\frac{\partial}{\partial x_{j}})+\sum_{i=1}^d
b_{i}\frac{\partial}{\partial x_{i}}
\end{equation}
with measurable coefficients $a:Q_T
\rightarrow{\mathbb{R}}^{d}\otimes{\mathbb{R}}^{d}$,
$b:Q_T\rightarrow{\mathbb{R}}^{d}$ such that
\begin{equation}\label{eq1.2}
\lambda|\xi|^2\leq\sum^d_{i,j=1}a_{ij}(t,x)\xi_i\xi_j\le\Lambda|\xi|^2,
\quad a_{ij}=a_{ji},
\quad|b_{i}(t,x)|\le \Lambda_{1}, \quad \xi\in{\mathbb{R}}^d
\end{equation}
for some  $0<\lambda\le\Lambda$ and $\Lambda_{1}>0$ (see
\cite{Y.Osh.DirichletForms,A.Roz.DIFF,Stroock.DIFF}).

It is known (see \cite{Le1,Roz.dec}) that for every $(s,x)\in
Q_{\hat{T}}\equiv [0,T)\times\mathbb{R}^{d}$ the  process
$X_{s,\cdot}\equiv X_{\cdot}-X_s$ is under $P_{s,x}$ a continuous
Dirichlet process on $[s,T]$ in the sense of F\"ollmer
\cite{Follmer}. In the paper we first develop some stochastic
calculus for time-dependent functionals of $\BX$. Secondly, we
give mild regularity conditions on $u$ under which the functional
$X^u_{s,\cdot}$ is a Dirichlet process under $P_{s,x}$ and, if it
is the case, we describe the martingale part $M^u$ and the
zero-quadratic variation part $A^u$ of its decomposition
\begin{equation}
\label{eq1.03} X^u_{s,t}=M^u_{s,t}+A^u_{s,t},\quad t\in[s,T],
\quad P_{s,x}\mbox{-}a.s.
\end{equation}
Finally, we characterize the class of $u$ such that
$X^u_{s,\cdot}$ is a semimartingale under $P_{s,x}$.

It is known that general Dirichlet processes are stable under
$C^{1}$ transformations (see \cite{Bertoin,CMS}). $C^1$-regularity
of $u$ is too strong in applications we have in mind. Our main
motivation to investigate functionals of the form $X^u$ comes from
the fact that they appear in probabilistic analysis of strong
solutions to parabolic PDEs or variational inequalities involving
the operator $L_t$ (see \cite{K,A.Roz.BSDE,Roz.Div}). Therefore
the natural assumption on $u$ is that it belongs to some Sobolev
space and in general is even not continuous.

Time-independent functionals of time-homogeneous diffiusions are
quite well investigated. Let $\XX$ be the locally compact
separable metric space and let $m$ a positive Radon measure on
$\XX$ such that supp$[m]=\XX$. Let $\{(X,P_{x});x\in\XX\}$ be an
$m$-symmetric Hunt process with Dirichlet form
$(\mathcal{E},D(\mathcal{E}))$ on $\mathbb{L}_{2}(\XX,m)$. It is
known (see \cite{Fukushima}) that for every $u\in D(\mathcal{E})$
there exists a $\mathcal{E}$-quasi-continuous version of $u$
(still denoted by $u$) such that $X^{u}$ admits the so called
Fukushima decomposition, i.e.
\[
X^{u}_{t}=M^{u}_{t}+A_{t}^{u},\quad t\in [0,T]\quad
P_{x}\mbox{-a.s.}
\]
for $\mathcal{E}-$q.e. $x\in\XX$, where $M^{u}$ is a continuous
martingale AF  of finite energy and $A^{u}$ is a continuous AF of
zero energy. An easy calculation (see \cite[pp. 201]{Fukushima})
shows that $A^{u}$ has zero-quadratic variation on $[0,T]$ under
the measure $P_{\nu}(\cdot)=\int_{\XX}P_x(\cdot)\,d\nu(x)$ along
dyadic partitions of $[0,T]$ for every Radon measure $\nu\ll m$.
Hence, to prove that $X^{u}$ is a Dirichlet process in the sense
of F\"ollmer one should relax the assumption on the absolute
continuity of $\nu$ and on the sequence of partitions. In
\cite{DMS} the authors weakened the assumption on the starting
measures $\nu$ in the case of Dirichlet form
$(\mathcal{E},D(\mathcal{E}))$ on
$\mathbb{L}_{2}(\mathbb{R}^{d},m)$ with the Lebegue measure $m$,
defined by
\begin{eqnarray}\label{D1}
\mathcal{E}(u,v)=\frac12\langle a\nabla u,\nabla
v\rangle_{2},\quad u, v\in D(\mathcal{E})=H^{1}(\mathbb{R}^{d}),
\end{eqnarray}
where $a(t,x)=a(x)$, $x\in\mathbb{R}^{d}$. The class of measures
considered in \cite{DMS} includes in particular the Dirac measure
$\delta_{\{x\}}$ for $\mathcal{E}$-q.e. $x\in\mathbb{R}^{d}$,
which shows that $X^{u}$ is a Dirichlet process on $[0,T]$ under
$P_{x}$ for $\mathcal{E}$-q.e. $x\in\mathbb{R}^{d}$ along dyadic
partitions. It is worth mentioning that in the case of
non-symmetric diffusions the approach of \cite{DMS} breaks down.

A different approach to the problem of investigating $X^u$ in case
\begin{equation}
\label{eq1.5} \mathcal{E}(u,v)=\frac12\langle a\nabla u,\nabla
v\rangle_{2} +\langle b\nabla u,v\rangle_{2}\quad u, v\in
D(\mathcal{E})=H^{1}(\mathbb{R}^{d}),
\end{equation}
with $a(t,x)=a(x)$, $b(t,x)=b(x)$ was adopted in
\cite{Roz.Dirichlet}. In \cite{Roz.Dirichlet} it is shown that if
$u\in W^{1}_{q}(\mathbb{R}^{d})$ with $q>2$ then $X^{u}$ is a
continuous Dirichlet process in the sense of F\"ollmer for
$\mathcal{E}$-q.e. $x\in\mathbb{R}^{d}$ (see also
\cite{A.Roz.DIFF,Roz.dec} where time-inhomogenous diffusions are
also considered).

In the case of one dimensional Wiener process $W$ it is known (see
\cite{FPS})  that $W^{u}$ is a continuous Dirichlet process in the
sense of F\"ollmer for every starting point $x\in\mathbb{R}$  if
$u\in H^{1}(\mathbb{R})$ and it appears that this condition is
necessary (see \cite{Chitashvili}). In the case of
multidimensional Wiener process one can deduce from \cite{FP} that
$W^{u}$ is a continuous Dirichlet process in the sense of
F\"ollmer on $[0,T]$ for q.e. starting points $x\in\mathbb{R}^{d}$
if $u\in H^{1}(\mathbb{R}^{d})$.

Up to our knowledge, in the case where $u$ depends on time, only
few results are available. In \cite{Le1} diffusions corresponding
to $L_t$ are considered. It is shown there that  $X^{u}$ is a
continuous Dirichlet process on $[s,T]$ in the sense of F\"ollmer
for every $(s,x)\in Q_{\hat{T}}$ if $\sup_{t\in[0,T]}(\|\nabla
u(t)\|_{p} +\|\frac{\partial u}{\partial t}(t)\|_{p})<\infty$ for
some $p>d\wedge 2$. In \cite{Chitashvili1} necessary and sufficient
conditions on $u$ for $X^u$ to be  semimartingale are given in
case $X$ is a one-dimensional Wiener process.

Let us now describe briefly the content of the paper. As already
mentioned, we are interested in solutions $u$ to parabolic PDEs or
parabolic variational inequalities involving $L_t$. Therefore our
basic assumption on $u$ is that $u\in{\mathcal{W}}_{\varrho}$,
where ${\mathcal{W}}_{\varrho}
=\{u\in{\mathbb{L}}_{2}(0,T;H^{1}_{\varrho}); \frac{\partial
u}{\partial t} \in{\mathbb{L}}_{2}(0,T;H^{-1}_{\varrho})\}$
($\varrho$ is some weight), i.e. $u$ belongs to the natural space
for strong solutions of such problems. Let
cap$_{L}:2^{Q_{\hat{T}}}\rightarrow \mathbb{R}^{+}\cup\{+\infty\}$
be the parabolic capacity associated with  $L_t$ (see
\cite{Pierre}) or, equivalently,  restriction to $Q_{\hat{T}}$ of
the capacity generated by time-dependent Dirichlet form
\[
\mathcal{E}(u,v)=\left\{
\begin{array}{l}
\int_{\mathbb{R}}\mathcal{E}^{(t)}(u(t),v(t))\,dt
-\int_{\mathbb{R}}\langle \frac{\partial u}{\partial t}(t),v(t)\rangle,\,\,
u\in\mathcal{W}, v\in\mathbb{L}_{2}(\mathbb{R};H^{1}),\medskip\\
\int_{\mathbb{R}}\mathcal{E}^{(t)}(u(t),v(t))\,dt
+\int_{\mathbb{R}}\langle \frac{\partial v}{\partial t}(t),u(t)\rangle,\,\,
v\in\mathcal{W}, u\in\mathbb{L}_{2}(\mathbb{R};H^{1}),
\end{array}
\right.
\]
where ${\mathcal{W}}$ denotes ${\mathcal{W}}_{\varrho}$ with
$\varrho\equiv1$,
\[
\mathcal{E}^{(t)}(u,v)=\frac12\langle a(t)\nabla u,\nabla
v\rangle_{2} +\langle b(t)\nabla u,v\rangle_{2},\quad t\in[0,T],
\]
$\mathcal{E}^{(t)}(u,v)=\mathcal{E}^{(0)}(u,v)$ for $t\le 0$ and
$\mathcal{E}^{(t)}(u,v)=\mathcal{E}^{(T)}(u,v)$ for $t\ge T$. In
the paper we provide various conditions on $u$ ensuring that for
cap$_{L}$-quasi every (q.e. for short) $(s,x)\in Q_{\hat{T}}$ the
process $X^{u}_{s,\cdot}$ is under $P_{s,x}$ a continuous
Dirichlet process on $[s,T]$ in the sense of F\"ollmer  or is a
continuous semimartingale.

For the convenience of the reader we begin in Section \ref{sec2}
with basic information on various definitions of parabolic
capacity associated with $L_t$.

In Section \ref{sec3} we formulate Fukushima's and the Lyons-Zheng
decomposition of $X$ under $P_{s,x}$. Using the latter
decomposition we investigate  additive functionals of the form
$\int\dyw\bar f(\theta,X_{\theta})\,d\theta$, where $\dyw\bar f$
stands for the divergence of the vector field $\bar
f=(f^1,\dots,f^d)$ such that $f^i\in\mathbb{L}_{2}^{loc}(Q_{T})$,
$i=1,\dots,d$. It is known that in case of time-homogeneous
diffusions $\{(X,P_x);x\in\mathbb{R}^d\}$ corresponding to $L_t$
with time-independent coefficients such functionals may be defined
under the measure $P_m$ as a forward-backward integral with
respect to martingales from the Lyons-Zheng decomposition of $X^u$
(see \cite{Stoica}). We show that the functionals can be well
defined for time-inhomogeneous diffusions and what is more
important, under the measure $P_{s,x}$ for q.e. $(s,x)\in Q_{\hat
T}$ (see \cite{Roz.Div} for similar results). We show also that if
$u\in H^{1}(\mathbb{R}^{d})$ then $X^{u}$ is a continuous
Dirichlet process in the sense of F\"ollmer under $P_{x}$ for
$\mathcal{E}$-q.e. $x\in\mathbb{R}^{d}$, where $\mathcal{E}$ is
given by (\ref{D1}).

In Section \ref{sec4} we show that each
$u\in{\mathcal{W}}_{\varrho}$ has a quasi-continuous version,
still denoted by $u$, such that $X^u_{s,\cdot}$ is a Dirichlet
process on $(s,T]$ under $P_{s,x}$ for every $(s,x)\in
Q_{\hat{T}}$. Under mild additional regularity conditions on $u$
it is a Dirichlet process on $[s,T]$ for cap$_{L}$-q.e. $(s,x)\in
Q_{\hat{T}}$. We describe also the martingale and the
zero-quadratic variation parts of the decomposition (\ref{eq1.03})
and show that (\ref{eq1.03}) implies the Fukushima decomposition
of $X^u$ into martingale AF of finite energy and CAF of zero
energy.

In Section \ref{sec5} we introduce the definition of the integral
with respect to continuous additive functionals (CAFs for short)
of $\BX$ of zero-quadratic variation associated with functionals
in ${\mathbb{L}}_{2}(0,T;H^{-1}_{\varrho})$. The key result here
says that given such CAF $A$ and bounded $\eta\in\WW_{\varrho}$
one can find a sequence $\{A^{n}\}$ of square-integrable CAFs of
finite variation such that for q.e. $(s,x)\in Q_{\hat{T}}$\,,
\[
E_{s,x}\sup_{s\le t \le T}|\int_{s}^{t}
\eta(\theta,X_{\theta})\,dA^{n}_{s,\theta}
-\int_{s}^{t}\eta(\theta,X_{\theta})\,dA_{s,\theta}|\rightarrow 0.
\]
This approximation result enables us to handle integrals with
respect to CAFs corresponding to functionals in
${\mathbb{L}}_{2}(0,T;H^{-1}_{\varrho})$. As a first application
we show that such CAFs are uniquely determined by their Laplace
transforms.

In Section \ref{sec6} we are concerned with the problem of finding
minimal conditions on $u\in{\mathcal{W}}_{\varrho}$ under which
$X^{u}$ is a semimartingale. Our main result proved here says that
$X^u_{s,\cdot}$ is a locally finite semimartingale under $P_{s,x}$
for q.e. $(s,x)\in Q_{\hat T}$ if and only if
$(\frac{\partial}{\partial t}+L_t)u$ is a signed Radon measure.

Finally, in Section \ref{sec7} we collect here some useful
estimates for diffusions ${\mathbb{X}}$ and related estimates on
the  fundamental solution $p$ and weak solutions of the Cauchy
problem associated with $L_t$.

In the paper we will use the following notation.

$Q_{st}=[s,t]\times{\mathbb{R}}^d$,
$Q_t=[0,t]\times{\mathbb{R}}^d$,
$Q_{\hat{T}}=[0,T)\times{\mathbb{R}}^d$,
$\nabla=(\frac{\partial}{\partial x_1},
\dots,\frac{\partial}{\partial x_d})$.
By $B(Q_{T})$ we denote the set of bounded Borel functions on
$Q_{T}$. $C_{c}(Q_{T})$ $(C_{c}({\mathbb{R}}^{d}))$ is the space
of all continuous functions with compact support in $Q_{T}$ (in
${\mathbb{R}}^{d}$).

${\mathbb{L}}_{p}({\mathbb{R}}^d)$ is the usual Banach space of
measurable functions  on ${\mathbb{R}}^d$ with the norm
$\|u\|_p=(\int_{{\mathbb{R}}^d}|u(x)|^p\,dx)^{1/p}$,
${\mathbb{L}}_{p,q}(Q_{tT})$ is the Banach space of measurable
functions on $Q_{tT}$ with the norm
$\|u\|_{p,q,t,T}=(\int_t^T(\int_{{\mathbb{R}}^d}
|u(s,x)|^p\,dx)^{p/q}\,ds)^{1/q}$,
${\mathbb{L}}_{p}(Q_{tT})={\mathbb{L}}_{p,p}(Q_{tT})$,
$\|u\|_{p,p,t,T}=\|u\|_{p,t,T}$ and $\|u\|_{p,T}=\|u\|_{p,0,T}$

Let $\varrho$ be a positive function on ${\mathbb{R}}^d$. By
${\mathbb{L}}_{p,\varrho}({\mathbb{R}}^d)$
(${\mathbb{L}}_{p,q,\varrho}(Q_{tT})$) we denote the space of
functions $u$ such that $u\varrho\in
{\mathbb{L}}_p({\mathbb{R}}^d)$
($u\varrho\in{\mathbb{L}}_{p,q}(Q_{t,T})$) equipped with the norm
$\|u\|_{p,\varrho}=\|u\varrho\|_p$
($\|u\|_{p,q,\varrho,t,T}=\|u\varrho\|_{p,q,t,T})$. We  write
$K\subset\subset X$ if $K$ is compact subset of $X$. By
$\langle\cdot,\cdot\rangle_2$ we denote the usual inner product in
${\mathbb{L}}_2({\mathbb{R}}^d)$ and by
$\langle\cdot,\cdot\rangle_{2,\varrho}$ the inner product in
${\mathbb{L}}_{2,\varrho}({\mathbb{R}}^d)$.

$W^{1}_{p,\varrho}$  is the Banach space consisting of all
elements $u$ of ${\mathbb{L}}_{p,\varrho}({\mathbb{R}}^d)$ having
generalized derivatives $\frac{\partial u}{\partial x_i}$,
$i=1,\dots,d$, in ${\mathbb{L}}_{p,\varrho}({\mathbb{R}}^d)$ with
dual  space $W^{-1}_{p,\varrho}$. We denote
$H^{1}_{\varrho}=W^{1}_{2,\varrho}$. ${\mathcal{W}}_\varrho$ is
the subspace of ${\mathbb{L}}_2 (0,T;H^1_{\varrho})$ consisting of
all elements $u$ such that $\frac{\partial u}{\partial t}
\in{\mathbb{L}}_2(0,T;H^{-1}_{\varrho})$, where $H^{-1}_{\varrho}$
is the dual space to $H^{1}_{\varrho}$ (see \cite{Lions} for
details). By  $\langle\cdot,\cdot\rangle_{\varrho}$ we denote the
duality pairing between spaces $H^{1}_{\varrho}$,
$H^{-1}_{\varrho}$ and by $\|\cdot\|_{*}$ we denote the norm in
Banach space ${\mathbb{L}}_2(0,T;H^{-1}_{\varrho})$.

By $\bar{f}$ we denote the vector function $(f^{1},\dots,f^{d})$.
We write $\bar f\in{\mathbb{L}}_{2,\varrho}(Q_{T})$ if $f^i\in
{\mathbb{L}}_{2,\varrho}(Q_{T})$, $i=1,\dots,d$. By $\mathcal{M}$
$({\mathcal{M}}^{+})$ we denote the set of all Radon measures
(positive Radon measures) on $Q_{T}$, and by
$\mathcal{M}^{+}([0,T])$ the set of positive  Radon measures on
$[0,T]$. $\mathcal{B}(E)$ $(\mathcal{B}_b(E),
\mathcal{B}^{loc}_{b}(E), \mathcal{B}^{+}(E))$ denotes the set of
all Borel (bounded, locally bounded, positive) real functions on a
topological space $E$.

By $C$ we  denote a general constant which may vary from line to
line but depends only on fixed parameters.

\nsubsection{Parabolic capacity} \label{sec2}

Let $\mathcal{R}$ denote the space of all measurable functions
$\varrho:\mathbb{R}^d\rightarrow\mathbb{R}$ such that
$\varrho(x)=(1+|x|^{2})^{-\alpha}$, $x\in\mathbb{R}^d$, for some
$\alpha\in\mathbb{R}$, and let $\mathcal{R}_{I}$ be the space of
all $\varrho\in\mathcal{R}$ such that
$\int_{\mathbb{R}^{d}}\varrho(x)\,dx<\infty$. Unless otherwise
stated, in the sequel we will always assume that
$\varrho\in\mathcal{R}_{I}$. We write also
$\varrho_{x}(y)=\varrho(y-x)$, $y\in{\mathbb{R}}^d$.

Let $\Phi\in{\mathbb{L}}_{2}(0,T;H^{-1}_{\varrho})$. It is well
known that $\Phi$ admit the decomposition
$\Phi=f^0+\mbox{div}\bar{f}$ for some $f^0,\bar{f}\in
{\mathbb{L}}_{2,\varrho}(Q_{T})$, i.e. $\Phi(\eta)=\langle
f^0,\eta\rangle_{2,\varrho}
-\langle\bar{f},\nabla(\varrho^{2}\eta)\rangle_{2}$. This
decomposition  is not unique but it is known that for every such
decomposition $\|\Phi\|_{*}\le
\|f^{0}\|_{2,\varrho,T}+\|\bar{f}\|_{2,\varrho,T}$ and there
exists a pair wich realizes the norm. If, in addition, $\Phi\ge
0$, i.e. $\Phi(\eta)\ge 0$ for any positive $\eta\in
{\mathbb{L}}_{2}(0,T;H^{1}_{\varrho})$, then by Riesz's theorem
there is a Radon measure $\mu$ on $Q_{T}$ such that
\begin{equation}\label{RM}
\Phi(\eta)=\int_{Q_{T}}\eta\,d\mu
\end{equation}
for every $\eta\in C_{0}^{\infty}(Q_{T})$. Let us observe that
\begin{eqnarray}\label{RM1}
\mu({\{t\}\times{\mathbb{R}}^{d}})=0, \quad t\in [0,T].
\end{eqnarray}
Indeed, if $\{\eta_{n}\}\subset
C_{0}^{\infty}(Q_{T})$ is a sequence of positive functions such
that $\eta_{n}\downarrow{\mathbf{1}}_{\{t\}\times{\mathbb{R}}^{d}}$
pointwise and in
${\mathbb{L}}_{2}(0,T;H^{1}_{\varrho})$, then
\[
0=\Phi(\eta)=\lim_{n\rightarrow\infty}\Phi(\eta_{n})
=\lim_{n\rightarrow\infty}\int_{\check{Q}_{T}}\eta_{n}\,d\mu
=\mu({\{t\}\times\mathbb{R}^{d}}).
\]

Let us define the capacity of $E\subset\subset
\check{Q}_{T}\equiv(0,T)\times{\mathbb{R}}^d$ by
\[
\overline{\mbox{cap}}_{\check{Q}_{T}}(E)
=\inf\{\int_{Q_{T}}|\nabla\eta(t,x)|^{2}\,dt\,dx: \eta\in
C_{0}^{\infty}(\check{Q}_{T}), \eta\geq {\mathbf{1}}_{E}\}.
\]
The capacity can be extended in a standard way to the Borel
$\sigma$-field ${\mathcal{B}}(\check{Q}_{T})$ of subsets of
$\check{Q}_{T}$. For $E\subset\subset \check{Q}_{T}$ and  $\eta\in
C^{\infty}_0(\check{Q}_T)$ such that $\eta \ge{\mathbf{1}}_E$ we
have
\[
\mu(E)\le \int_{{Q}_{T}}\eta \,d\mu =\Phi(\eta)=\int_{Q_{T}}\eta
f^0\varrho^{2}-\int_{Q_{T}}\bar{f}\nabla (\varrho^2 \eta)\le
C\|\nabla \eta\|_{2,\varrho,T}\,,
\]
the last inequality being a consequence of the
Gagliardo-Nirenberg-Sobolev inequality (see, e.g., \cite{Evans1}).
Thus, $\mu\ll\overline{\mbox{cap}}_{\check{Q}_{T}}$. Now, for
$E\subset\subset {\mathbb{R}}^{d}$ define
\[
\mbox{cap}_{{\mathbb{R}}^{d}}(E)=\inf\{\int_{{\mathbb{R}}^{d}}
|\nabla\eta(x)|^{2}\,dx:\eta\in C_{0}^{\infty}({\mathbb{R}}^{d}),
\eta\geq{\mathbf{1}}_{E}\},
\]
and extend it in the standard way to
${\mathcal{B}}({\mathbb{R}}^{d})$. From \cite{B.Mo.} it follows
that for every $B\in{\mathcal{B}}(\check{Q}_{T})$,
\[
\overline{\mbox{cap}}_{\,\check{Q}_{T}}(B)=\int_{0}^{T}
\mbox{cap}_{\,{\mathbb{R}}^d}(B_{t})\,dt,
\]
where  $B_{t}=\{x\in{\mathbb{R}}^{d};(t,x)\in B\}$. Since
$\mu\ll\overline{\mbox{cap}}_{\check{Q}_{T}}$, using the the well
known fact that elements of $H^{1}_{\varrho}$ have
quasi-continuous versions defined up to the sets of
$\mbox{cap}_{{\mathbb{R}}^{d}}$-measure zero (see \cite[Chapter
2]{Fukushima} we may extend formula (\ref{RM}) to  all $\eta\in
{\mathbb{L}}_{2}(0,T;H^{1}_{\varrho})$.

It is worth noting that in the definition of capacity
$\overline{\mbox{cap}}_{Q_{T}}$ and in the representation theorem
for functionals in ${\mathbb{L}}_{2}(0,T;H^{-1}_{\varrho})$
derivatives with respect to the time variable do not appear.
Therefore various facts on functionals $\mu\in
{\mathbb{L}}_{2}(0,T;H^{-1}_{\varrho})\cap{\mathcal{M}}$ can be
proved by making obvious changes in proofs of corresponding facts
concerning elliptic capacity and  functionals in
$H^{-1}_{\varrho}$. In particular, modifying slightly arguments
from \cite{Boccardo} and \cite{DalMaso} one can prove the following theorems.

\begin{theorem}\label{decom}
Let $\mu\in \mathcal{M}$. If $\mu\ll\overline{\mbox{\rm
cap}}_{Q_{T}}$ then there exists $\gamma_{1},\gamma_{2}\in
{\mathbb{L}}_{2}(0,T;H^{-1}_{\varrho})\cap{\mathcal{M}}^{+}$ and
positive $\alpha_{i}\in{\mathbb{L}}_{1,loc}(Q_{T},\gamma_{i})$,
$i=1,2$ such that
$d\mu=\alpha_{1}\,d\gamma_{1}-\alpha_{2}\,d\gamma_{2}$.
\end{theorem}

\begin{theorem}\label{decom1}
A Radon measure $\nu$ vanishes on sets of zero $\overline{{\rm
cap}}_{Q_{T}}$  capacity if and only if it admits a decomposition
\[
\mu=\Phi+k,
\]
where $\Phi\in \mathbb{L}_{2}(0,T;H^{-1}_{\varrho})$  and $k\in
\mathbb{L}_{1}^{loc}(Q_{T})$.
\end{theorem}

In the paper we will use also another notion of capacity, the so
called parabolic capacity, which appears when considering the
natural space of strong solutions of variational inequalities,
i.e. the space ${\mathcal{W}}_{\varrho}$.

Let $\Omega=C([0,T],{\mathbb{R}}^d)$ denote the space of
continuous ${\mathbb{R}}^d$-valued functions on $[0,T]$ equipped
with the topology of uniform convergence and let $X$ be the
canonical process on $\Omega$. It is known that for given operator
$L_t$ defined by (\ref{eq1.1}) with $a$ and $b$ satisfying
(\ref{eq1.2}) one can construct a weak fundamental solution $p$
for $L_t$ and then a Markov family
${\mathbb{X}}=\{(X,P_{s,x});(s,x)\in Q_{\hat T}\}$ for which $p$
is the transition density function, i.e.
\[
P_{s,x}(X_t=x;0\leq t\leq s)=1,\quad
P_{s,x}(X_t\in\Gamma)=\int_{\Gamma}p(s,x,t,y)\,dy,\quad t\in(s,T]
\]
for any $\Gamma\in{\mathcal{B}}({\mathbb{R}}^d)$ (see
\cite{A.Roz.DIFF,Stroock.DIFF}). We define the parabolic capacity
of a Borel set $B\subset Q_{\hat{T}}$ by
\[
\mbox{cap}_{L}(B)=P_{m}(\exists\,t\in [s,T):(t,X_{t})\in B),
\]
where $m$ is the Lebesgue measure on ${\mathbb{R}}^d$ and
\[
P_{m}(\Gamma)=\int_{Q_{\hat{T}}}P_{s,x}(\Gamma)\,ds\,dx,\quad
\Gamma\in{\mathcal{B}}({\mathbb{R}}^d).
\]

In what follows we say that some property is satisfied
quasi-everywhere (q.e. for short) if it is satisfied except of a
Borel set of zero capacity $\mbox{cap}_{L}$.

\begin{remark}\label{rem.76}
{\rm It follows directly from the definition of 
$\mbox{cap}_{L}$ that cap$_{L}(\{s\}\times B)>0$ for every $s\in
(0,T)$ and $B\in\mathcal{B}(Q_{T})$ such that $m(B)>0$. Hence, if
some property holds for q.e. $(s,x)\in Q_{\hat{T}}$, then it holds
for a.e. $x\in\mathbb{R}^{d}$ for every $s\in (0,T)$. }
\end{remark}

From \cite{Oshima1} and \cite{Pierre} it follows that the
parabolic capacity $\mbox{cap}_{L}$ is equivalent to the following
parabolic capacity $\mbox{cap}_{2}$ in the analytical sense.

\begin{definition}
{\em Let $V\subset Q_{\hat{T}}$ be an open set. We set
\[
\mbox{cap}_{2}(V)=\inf\{\|u\|_{\WW_{\varrho}}: u \in\WW_{\varrho},
u\ge {\mathbf{ 1}}_{V}\mbox{ a.e.}\}
\]
with the convention that $\inf\emptyset=\infty$. The parabolic
capacity of a Borel $B\subset Q_{\hat{T}}$ is defined by
\[
\mbox{cap}_{2}(B)=\inf\{\mbox{cap}_{2}(V): V \mbox{ is an open
subset of }Q_{T}, B\subset V\}.
\]
}
\end{definition}

From \cite[Proposition 2]{Pierre} it follows that $\mbox{cap}_{2}$
is a Choquet capacity. In particular, it follows (see
\cite[Theorem A.1.1.]{Fukushima}) that for any
$B\in\mathcal{B}(Q_{T})$,
\begin{equation}\label{eq.cap}
\mbox{cap}_{2}(B)=\sup_{K\subset B,
K-\mbox{compact}}\mbox{cap}_{2}(K).
\end{equation}

\begin{remark}
\label{choquet} {\rm Since $\mbox{cap}_{2}$ and $\mbox{cap}_{L}$
are equivalent, it follows from (\ref{eq.cap}) that
$\mbox{cap}_{L}(B)=0$ if $\mbox{cap}_{L}(K)=0$ for every compact
subset $K$ of $B$.}
\end{remark}

\begin{definition}
{\rm We say that $u:Q_{T}\rightarrow \mathbb{R}$ is
quasi-continuous if $u$ is Borel measurable and $[0,T]\ni t\mapsto
u(t,X_{t})$ is a continuous process under the measure $P_{s,x}$
for q.e. $(s,x)\in Q_{\hat{T}}$\,. }
\end{definition}

The  notion of quasi-continuity defined above is equivalent to the
following one: for every every $\varepsilon>0$ there exists an
open set $U_{\varepsilon}\subset\check{Q}_{T}$ such that
$u_{|\check{Q}_{T}\setminus U_{\varepsilon}}$ is continuous  and
cap$_{2}(U_{\varepsilon})<\varepsilon$ (see Remark 3.6 and
Proposition 3.7 in \cite{Stannat1}). Let us note also that it is
known that every $u\in\mathcal{W}_{\varrho}$ has a
quasi-continuous version (see \cite{Oshima1}).

\nsubsection{Diffusions corresponding to divergence form
operators} \label{sec3}

Set ${\mathcal{F}}^s_t=\sigma(X_u,u\in[s,t]),
\bar{\mathcal{F}}^{s}_{t}=\sigma(X_u,u\in[T+s-t,T])$ and define
${\mathcal{G}}$ as the completion of ${\mathcal{F}}^s_T$ with
respect to the family ${\mathcal{P}}=\{P_{s,\mu}:\mu$ is a
probability measure on ${\mathcal{B}}({\mathbb{R}}^d)$\}, where
$P_{s,\mu}(\cdot)=\int_{{\mathbb{R}}^d}P_{s,x}(\cdot)\,\mu(dx)$,
and define ${\mathcal{G}}^s_t$ ($\bar{\mathcal{G}}^s_t)$ as the
completion of ${\mathcal{F}}^s_t$ ($\bar{\mathcal{F}}^s_t)$ in
${\mathcal{G}}$ with respect to ${\mathcal{P}}$.

We will say that a family $A=\{A_{s,t}, 0\le s\le t\le T\}$ of
random variables is an additive functional (AF) of ${\mathbb{X}}$
if $A_{s,\cdot}$ is $\{{\mathcal{G}}_{t}^{s}\}$-measurable for
every $0\le s\le t\le T$ and $P_{s,x}(A_{s,t}=A_{s,u}+A_{u,t},
s\le u\le t\le T)=1$ for q.e. $(s,x)\in Q_{\hat{T}}$. If, in
addition, $A_{s,\cdot}$ has $P_{s,x}$-almost all continuous
trajectories for q.e. $(s,x)\in Q_{\hat{T}}$, then  $A$ is called
a continuous AF (CAF), and if  $A_{s,\cdot}$ is an increasing
process under $P_{s,x}$ for  q.e. $(s,x)\in Q_{\hat{T}}$, it  is
called an increasing AF or positive AF. If $M$ is an AF such that
for q.e. $(s,x)\in Q_{\hat{T}}$, $E_{s,x}|M_{s,t}|^2<\infty$ and
$E_{s,x}M_{s,t}=0$ for $t\in[s,T]$ ($E_{s,x}$ is the expectation
with respect to $P_{s,x}$), it is called a martingale AF (MAF). We
say that $A$ is an  AF (CAF, increasing AF, MAF) in the strict
sense if the corresponding property holds for every $(s,x)\in
Q_{\hat{T}}$. Finally, we say that $A$ is a quasi-strict AF (CAF,
increasing AF, MAF) if the corresponding property holds under
$P_{s,x}$ for every $(s,x)\in Q_{\hat{T}}$ on $(s,T]$ and for q.e.
$(s,x)\in Q_{\hat T}$ on $[s,T]$. Since in  what follows, except
of Proposition \ref{propl.1},  we will consider exclusively
quasi-strict AFs, we will call it briefly additive functionals.

\subsubsection{Fukushima's decomposition and decomposition in
the sense of F\"ollmer}

It is known (see \cite{Le1,Roz.dec}) that there exist CAF $A$ in
the strict sense and a continuous MAF $M$ in the strict sense such
that
\begin{equation}\label{eq3.1}
X_t-X_s=M_{s,t}+A_{s,t},\quad t\in[s,T],\quad P_{s,x}\mbox{-}a.s.,
\end{equation}
for every $(s,x)\in Q_{\hat T}$, and moreover, $M_{s,\cdot}$ is
a $(\{{\mathcal{G}}^s_t\},P_{s,x})$-square-integrable martingale
on $[s,T]$ with the co-variation  given by
\begin{equation}
\label{eq2.3} \langle M^{i}_{s,\cdot},M^{j}_{s,\cdot}\rangle_t=
\int_s^ta_{ij}(\theta, X_\theta)\,d\theta,\quad t\in[s,T],\quad
i,j=1,...,d,
\end{equation}
while $A_{s,\cdot}$ is a process of $P_{s,x}$-zero-quadratic
variation on $[s,T]$, i.e. $A_{s,s}=0$ and
\begin{eqnarray}\label{eqp.27}
\sum_{t_i\in{\mit\Pi}_m} |A_{s,t_{i+1}}-A_{s,t_i}|^2\rightarrow0
\mbox{ in probability $P_{s,x}$}
\end{eqnarray}
for any sequence $\{{\mit\Pi}_m=\{t_0,t_1,\dots,t_{i(m)}\}\}$ of
partitions of $[s,T]$ such that $s=t_0<t_1<\dots<t_{i(m)}=T$ and
$\|\mit\Pi_m\|=\max_{1\le i\le i(m)}|t_i-t_{i-1}|\rightarrow$ $0$
as $m\rightarrow\infty$. In particular, $X_{\cdot}-X_s$ is a
$(\{{\mathcal{G}}^s_t\},P_{s,x})$-Dirichlet process in the sense
of F\"ollmer. One can show also that $M$ is a MAF of locally
zero-energy and $A$ is a CAF of locally finite energy (see
\cite{Roz.dec} and \cite{R,Roz.Dirichlet} for time-homogeneous
diffusions), i.e. (\ref{eq3.1}) coincides with Fukushima's
decomposition for ${\mathbb{X}}$.

Observe that if $\sigma\sigma^*=a$ then by (\ref{eq2.3}),
\begin{equation}
\label{eq2.5}
B_{s,t}=\int^t_s\sigma^{-1}(\theta,X_{\theta})\,dM_{s,\theta},\quad
t\in[s,T]
\end{equation}
is a $(\{{\mathcal{G}}^s_t\},P_{s,x})$-Wiener process.

\subsubsection{The Lyons-Zheng decomposition}

Additional information on the structure of $A$ of decomposition
(\ref{eq3.1}) provides the Lyons-Zheng decomposition for
${\mathbb{X}}$. Let $(s,x)\in Q_{\hat{T}}$. For $s\le u\le t \le
T$ we set
\[
\alpha^{s,x,i}_{u,t}=\sum_{j=1}^d
\int_{u}^{t}\frac{1}{2}a_{ij}(\theta,X_{\theta})p^{-1}
\frac{\partial p}{\partial y_j}
(s,x,\theta,X_{\theta})\,d\theta,\quad
\beta^{i}_{u,t}=\int_{u}^{t}b^{i}(\theta,X_{\theta})\,d\theta.
\]

In the sequel, for a process $Y$ on $[s,T]$ and fixed measure
$P_{s,x}$ we  write $\bar{Y}_t=Y_{T+s-t}$ for $t\in[s,T]$.

From \cite{Roz.dec} it follows that under $P_{s,x}$ the canonical
process $X$ admits the decomposition
\begin{equation}\label{dec1}
X_{t}-X_{u}=\frac{1}{2}M_{u,t}+\frac{1}{2}
(N^{s,x}_{s,T+s-t}-N^{s,x}_{s,T+s-u})
-\alpha^{s,x}_{u,t}+\beta_{u,t},\quad s\le u\le t\le T,
\end{equation}
where  $M_{s,\cdot}$ is the martingale of (\ref{eq3.1}) and
$N^{s,x}_{s,\cdot}$ is a
$(\{\bar{\mathcal{G}}^{s}_{t}\},P_{s,x})$-martingale such that
\begin{equation}
\label{eq2.05} \langle
N^{s,x,i}_{s,\cdot},N^{s,x,j}_{s,\cdot}\rangle_t = \int_s^t
a_{ij}(\bar{\theta}, \bar{X}_\theta)\,d\theta,\quad
t\in[s,T],\quad i,j=1,\dots,d.
\end{equation}
Observe that co-variation of $N^{s,x}$ does not depend on
$x\in{\mathbb{R}}^d$.

\begin{remark}
{\rm From (2.7) in \cite{A.Roz.DIFF} it follows that
\[
E_{s,x}(\tilde
N^{s,x}_{u_1,u_2}|\bar{\mathcal{G}}^{s}_{T+s-u_2})=0,\quad s\le
u_1\le u_2\le T.
\]
Hence, if we put
\[
\bar{M}_{u,t}^{s,x}=-(N^{s,x}_{s,T+s-t}-N^{s,x}_{s,T+s-u}),
\]
then for every $t\in[s,T)$, $\{\bar{M}_{t+s-u,t},u\in[s,t]\}$ is a
$(\{\bar{\mathcal{G}}^{t+s-u}_{T}\}_{u\in[s,t]},P_{s,x})$-martingale
and under $P_{s,x}$ the process $X$ admits the decomposition
\[
X_{t}-X_{u}=\frac12 M_{u,t}-\frac12\bar{M}^{s,x}_{u,t}
-\alpha^{s,x}_{u,t}+\beta_{u,t},\quad s\le u\le t\le T
\]
considered in \cite{Lyons}. }
\end{remark}

\subsubsection{Forward-backward integrals}

Let $\bar f=(f_1,\dots,f_d):Q_T\rightarrow\BR^d$ and let $S$ be
some class of real functions defined on $Q_T$. To simplify notation, in
what follows we write $\bar{f}\in S$ if $f_i\in S$, $i=1,\dots,d$.

Let $\bar{f}\in\mathcal{B}^{loc}_{b}(Q_{T})$. Similarly to
\cite{Roz.Div,Stoica}, using (\ref{dec1}) we set under the measure
$P_{s,x}$\,,
\begin{eqnarray}
\label{eq2.6}
\int_{r}^{t}\bar{f}(\theta,X_{\theta})\,d^{*}X_{\theta}\equiv
-\int_{r}^{t}\bar{f}
(\theta,X_{\theta})(dM_{s,\theta}+d\alpha^{s,x}_{s,\theta})
-\int^{T+s-r}_{T+s-t}
\bar{f}(\bar{\theta},\bar{X}_{\theta})dN^{s,x}_{s,\theta}
\end{eqnarray}
for $s\le u\le t\le T$, where $\bar\theta=T+s-\theta$. By
Proposition \ref{prop3.2}, all integrals on the right-hand side of
(\ref{eq2.6})  are well defined for every $(s,x)\in Q_{\hat{T}}$.
The interest in the integral defined above comes from the fact
that if $\bar f$ is regular then
\begin{equation}
\label{eq2.7}
\int_{u}^{t}\mbox{div}\bar{f}(\theta,X_{\theta})\,d\theta
=\int_{u}^{t}a^{-1}\bar{f}(\theta,X_{\theta})\,d^{*}X_{\theta},\quad
s\le u\le t\le T,\quad P_{s,x}\mbox{-}a.s.
\end{equation}
(see \cite{Roz.Div}), which enables one  to extend the integral on
the left-hand side of (\ref{eq2.7}) to
$\bar{f}\in\mathcal{B}^{loc}_b(Q_{T})$.

Our first goal is to extend the class of functions for which
(\ref{eq2.6}) is well defined for q.e. $(s,x)\in Q_{\hat T}$. In
view of (\ref{eq1.2}), (\ref{eq2.3}), (\ref{eq2.05})  and
Proposition \ref{prop3.2}, to define integrals with respect to the
forward and backward martingales  it suffices to assume that
$f\in\mathbb{L}_{2,\varrho}(Q_{T})$. The main problem is to define
integral with respect to $\alpha^{s,x}$, because the gradient of
$p$ is not square-integrable (see \cite{Aronson}) and
$\alpha^{s,x}$ depends on $(s,x)$. The latter fact makes
difficulties  in applying the Markov property of $\BX$ to get
existence of the integral.

We start with the investigation of integrals with respect to
$\alpha^{s,x}$ in case of time-homogeneous diffusions. Let
$\{(X,P_{x});x\in\mathbb{R}^d\}$ be a Hunt process associated with
the Dirichlet form (\ref{eq1.5}).

It is known that if $a$ is piecewise smooth (see \cite{DKN} for
details) then  there exists $M>0$ such that
\[
|\nabla_{x} p(t,x,y)|\le
\frac{M}{t^{(d+1)/2}}\exp(-\frac{|y-x|^{2}}{2Mt}).
\]
Hence, by the elementary calculations,
\[
E_{x}\int_{0}^{T}|f(X_t)|\,d|\alpha^{x}|_{t} \le
C\int_{\mathbb{R}^{d}}f(y)|y-x|^{1-d}\,dy.
\]
On the right-hand side of the above inequality we recognize the
Riesz potential of order $1$.  Therefore repeating arguments from
the  proof of \cite[Proposition 3.6]{FP} shows that for every
$f\in{\mathbb{L}}_{2,loc}(\mathbb{R}^{d})$,
\begin{equation}\label{eql.1}
P_{x}(\int_{0}^{T}|f(X_{t})|\,d|\alpha^{x}|_{t}<\infty)=1
\end{equation}
for $\mathcal{E}$-q.e. $x\in\mathbb{R}^{d}$.

The following example shows that in the time-dependent case  the
condition $f\in\mathbb{L}_{2,\varrho}(Q_{T})$ is insufficient to
guarantee (\ref{eql.1}) even if $a$ is smooth.

\begin{example}
{\rm Let $d=1$, $a=1,b=0$, so
that $X_t$, $t>s$, has  under $P_{s,x}$ the normal distribution
with mean $x$ and variance  $t-s$. Then
\[
\alpha^{s,x}_{s,t}=\frac12\int_{s}^{t}p^{-1}\frac{\partial
p}{\partial y} (s,x,\theta,X_{\theta})\,d\theta
=\frac12\int^t_s\frac{x-X_{\theta}}{\theta-s}\,d\theta.
\]
Suppose that $f$ is nonnegative and does not depend on $x$. Then
\[
w(s,x)=E_{s,x}\int_{s}^{T}f(\theta)\,d|\alpha^{s,x}_{s,\cdot}|_{\theta}
= E_{s,x}\int_{s}^{T}f(\theta)\frac{|x-X_{\theta}|}{\theta-s}\,d\theta
=C\int_{s}^{T}\frac{f(\theta)}{(\theta-s)^{1/2}}\,d\theta,
\]
i.e. $w(s,x)$ does not depend on $x$. Now, let us fix $t_{0}\in
(0,T)$. Since the function $(t_0,T)\ni t\mapsto(t-t_0)^{-1/2}$
does not belong to ${\mathbb{L}}_{2}(t_{0},T)$, one can find
$f\in{\mathbb{L}}_{2}(0,T)$ such that
$\int_{t_{0}}^{T}f(\theta)(\theta-t_{0})^{-1/2}=\infty$.  Then
$w=\infty$ on the set $\{t_{0}\}\times\mathbb{R}^{d}$ and  from
Remark \ref{rem.76} it follows that
cap$_{L}(\{t_{0}\}\times\mathbb{R}^{d})>0$. }
\end{example}

We will extend the integral side of (\ref{eq2.6}) to $\bar{f}\in
{\mathbb{L}}_{2,\varrho}(Q_{T})$ by using approximation.

\begin{proposition}\label{propl.1}
Let  $p>0$ and let $A,A^{n}$, $n\in\mathbb{N}$, be CAFs of $\BX$
such that
\begin{eqnarray}\label{eql.2}
E_{s,x}\sup_{s\le t\le T}|A^{n}_{s,t}-A_{s,t}|^{p}\rightarrow 0
\end{eqnarray}
for a.e. $(s,x)\in Q_{\hat T}$. Then there exists a subsequence
$\{n'\}$ of $\{n\}$ such (\ref{eql.2}) holds along $\{n'\}$ for
q.e. $(s,x)\in Q_{\hat T}$.
\end{proposition}
\begin{proof}
Set $B=\{(s,x):E_{s,x}\sup_{s\le t\le T}|
A^{n}_{s,t}-A_{s,t}|^{p}\nrightarrow 0\}$ and let $\tau=\inf
\{t\in[s,T]: (t,X_{t})\in K\}$, where $K$ is a compact subset of
$B$. Since $(X,P_{s,x})$ is a Feller process, $\tau$ is a
$\{{\mathcal{G}}^{s}_{t}\}$-stopping time. Hence, by the strong
Markow property with random shift and additivity of $A^{n}$ and
$A$,
\begin{align*}
P_{s,x}(\tau<\infty)&=P_{s,x}(E_{\tau\wedge T,X_{\tau\wedge T}}
\sup_{\tau\wedge T\le t\le T}| A^{n}_{\tau\wedge T,t}-A_{\tau\wedge T,t}|^{p}
\nrightarrow 0)\\
&=P_{s,x}(E_{s,x}(\sup_{\tau\wedge T\le t\le T}| A^{n}_{\tau\wedge T,t}
-A_{\tau\wedge T,t}|^{p}|{\mathcal{G}}^{s}_{\tau\wedge T})
\nrightarrow 0)\\& \le P_{s,x}(2^{p\wedge 1}
E_{s,x}(\sup_{s\le t\le T}| A^{n}_{s,t}
-A_{s,t}|^{p}|{\mathcal{G}}^{s}_{\tau\wedge T})\nrightarrow 0).
\end{align*}
Set
\[
T_{n,m}(s,x,\omega)=E_{s,x}(\sup_{s\le t\le T} |
A^{n}_{s,t}-A_{s,t}|^{p}|{\mathcal{G}}^{s}_{\tau\wedge
T})(\omega).
\]
By (\ref{eql.2}), $T_{n,m}\rightarrow 0$ in
${\mathbb{L}}_{1}(Q_{T}\times\Omega,\Pi)$, where $\Pi$ is the
finite measure defined by the formula
\[
\Pi(B)=\int_{Q_{T}}(E_{s,x}{\mathbf{1}}_{B}(s,x))\varrho(x)\,ds\,dx.
\]
Using the Borel-Cantelli lemma we can chose a subsequence (still
denoted $n$)  such that $T_{n,m}\rightarrow 0$, $\Pi$-a.e.. In
particular, $T_{n,m}(s,x)\rightarrow 0$, $P_{s,x}$-a.s. for a.e.
$(s,x)\in Q_{T}$. Hence $P_{s,x}(\tau<\infty)=0$ for a.e.
$(s,x)\in Q_{T}$, and consequently $\mbox{cap}_L(K)=0$. Hence, by
Remark \ref{choquet}, $\mbox{cap}_{L}(B)=0$.
\end{proof}

\begin{corollary}\label{coll.1}
Let  $p>0$ and let $A$ be a CAF of $\BX$ such that
\begin{equation}\label{eql.3}
E_{s,x}\sup_{s\le t\le T}|A_{s,t}|^{p}<\infty
\end{equation}
for a.e. $(s,x)\in Q_{\hat T}$. Then (\ref{eql.3}) holds for q.e.
$(s,x)\in Q_{\hat T}$.
\end{corollary}

\begin{proposition}
\label{prop.cd} Let $\bar{f}\in{\mathbb{L}}_{2,\varrho}(Q_{T})$.
Then there exists a unique   CAF $D$ of $\BX$
such that for every sequence $\{\bar
f_n\}\subset\mathcal{B}_b(Q_{T})$ convergent to $\bar{f}$ in
${\mathbb{L}}_{2,\varrho}(Q_{T})$ there exists a subsequence
(still denoted $n$) such that
\begin{equation}
\label{eql.4} E_{s,x}\sup_{s\le t\le T}|\int_{s}^{t}
\bar{f}_{n}(\theta,X_{\theta})\,d^{*}X_{\theta}-D_{s,t}|\rightarrow
0
\end{equation}
for q.e. $(s,x)\in Q_{\hat{T}}$.
\end{proposition}
\begin{proof}
Put
$A^n_{s,t}=\int_{s}^{t}\bar{f}_{n}(\theta,X_{\theta})\,d^{*}X_{\theta}$.
By Proposition \ref{prop3.2},
\begin{equation}
\label{eq2.12} \int_{Q_{T}}(E_{s,x}\sup_{s\le t \le T}|
A^{n}_{s,t} -A^{m}_{s,t}|)\varrho(x)\,dx \rightarrow 0.
\end{equation}
Hence, by Proposition \ref{propl.1}, there exist a subsequence
(still denoted $\{n\}$) and some process $D^{s,x}$ such that
$E_{s,x}\sup_{s\le t \le T}| A^{n}_{s,t} -D^{s,x}_{s,t}|
\rightarrow 0$ for q.e. $(s,x)\in Q_{\hat{T}}$. Using arguments
from the proof of \cite[Lemma A.3.2]{Fukushima} one can choose a
version of $D^{s,x}$ which does not depend on $(s,x)$.

To prove uniqueness, suppose that $\tilde{D}$ is another process
having the properties of $D$. Let $\bar{g}_{n}, \bar{f}_{n}\in
\mathcal{B}_{b}(Q_{T})$, $\bar{f}_{n},\bar{g}_{n}\rightarrow
\bar{f},\bar{g}$ in $\mathbb{L}_{2,\varrho}(Q_{T})$. Let $\{n\}$
be a subsequence such that (\ref{eql.4}) holds with the pairs
$(\bar{f}_{n},D)$, $(\bar{g}_{n},\tilde{D})$ and
$(\bar{f}_{n},\bar{g}_{n})$. For the latter pair it is possible
thanks to Proposition \ref{propl.1} and the the following
convergence
\[
\int_{Q_{T}}(E_{s,x}\sup_{s\le t \le T}|
\int_{s}^{t}\bar{f}_{n}(\theta,X_{\theta})\,d^{*}X_{\theta}
-\int_{s}^{t}\bar{g}_{n}(\theta,X_{\theta})\,d^{*}X_{\theta}|)\varrho(x)\,dx
\rightarrow 0,
\]
which is a consequence of convergence of $\{\bar{f}_{n}\},
\{\bar{g}_{n}\}$  and Proposition \ref{prop3.2}. Finally for q.e.
$(s,x)\in Q_{\hat{T}}$
\begin{align*}
E_{s,x}\sup_{s\le t\le T}|D_{s,t}-\tilde{D}_{s,t}| & \le
\lim_{n\rightarrow\infty} E_{s,x}\sup_{s\le t\le T} |\int_{s}^{t}
\bar{f}_{n}(\theta,X_{\theta})\,d^{*}X_{\theta}-D_{s,t}|\\
&\quad +\lim_{n\rightarrow\infty} E_{s,x}\sup_{s\le t \le T}|
\int_{s}^{t}\bar{f}_{n}(\theta,X_{\theta})\,d^{*}X_{\theta}
-\int_{s}^{t}\bar{g}_{n}(\theta,X_{\theta})\,d^{*}X_{\theta}| \\
&\quad +\lim_{n\rightarrow\infty} E_{s,x}\sup_{s\le t\le
T}|\int_{s}^{t}
\bar{g}_{n}(\theta,X_{\theta})\,d^{*}X_{\theta}-\tilde{D}_{s,t}|=0,
\end{align*}
and the proof is complete.
\end{proof}

\begin{remark}\label{reml.1}
{\rm For every $(s,x)\in Q_{\hat{T}}$ and $s< r\le t\le T $ the
integrals on the right-hand side of (\ref{eq2.6}) are well defined
$P_{s,x}$-a.s. This follows from Aronson's estimates and
Proposition \ref{prop.p}(ii), because
\begin{align}\label{inl.1}
E_{s,x}\int_{r}^{T}|\bar{f}(\theta,X_{\theta})|^2\,d\theta
&=\int_{Q_{rt}}|\bar{f}(\theta,y)|^2
p(s,x,\theta,y)\,d\theta\,dy\nonumber\\
&\le C\int_{Q_{rT}}\frac{1}{(\theta-s)^{d/2}}|\bar{f}(\theta,y|^2
\exp(\frac{-|y-x|}{C(\theta-s)})\nonumber\\
&\le C\frac{\varrho^{-1}(x)}{(r-s)^{d/2}}\|\bar{f}\|_{2,\varrho,T}
\end{align}
and
\begin{align}\label{inl.2}
E_{s,x}\int_{r}^{t}
|\bar{f}(\theta,X_{\theta})|\,d|\alpha^{s,x}_{s,\cdot}|_{\theta}&\le
\int_{Q_{r,T}}|\bar{f}(\theta,y)
|\nabla_{x}p|(s,x,\theta,y)|\,d\theta\,dy\nonumber\\
&\le \|\bar{f}\|_{2,\varrho,T}\,\|\nabla
p(s,x)\|_{2,\varrho^{-1},r,T}\,.
\end{align}
}
\end{remark}

\begin{proposition}\label{lml.1}
Let $\bar{f}\in \mathbb{L}_{2,\varrho}(Q_{T})$ and let $D$ be the
CAF of Proposition \ref{prop.cd}. Then $P_{s,x}$-a.s.,
\begin{eqnarray*}
D_{r,t}=-\int_{r}^{t}\bar{f}
(\theta,X_{\theta})(dM_{s,\theta}+d\alpha^{s,x}_{s,\theta})
-\int^{T+s-r}_{T+s-t}
\bar{f}(\bar{\theta},\bar{X}_{\theta})\,dN^{s,x}_{s,\theta},\quad
s<r\le t\le T
\end{eqnarray*}
for q.e. $(s,x)\in Q_{\hat{T}}$.
\end{proposition}
\begin{proof}
Let $\{\bar{f}_{n}\}\subset \mathcal{B}_b(Q_{T})$ be such  that
(\ref{eql.4}) holds q.e.. Then by (\ref{eq2.6}),
\begin{eqnarray*}
\int_{r}^{t}\bar{f}_{n}(\theta,X_{\theta})\,d^{*}X_{\theta}=
-\int_{r}^{t}\bar{f}_{n}
(\theta,X_{\theta})(dM_{s,\theta}+d\alpha^{s,x}_{s,\theta})
-\int^{T+s-r}_{T+s-t}
\bar{f}_{n}(\bar{\theta},\bar{X}_{\theta})\,dN^{s,x}_{s,\theta},
\end{eqnarray*}
and the result follows from (\ref{eql.4}), (\ref{inl.1}) and
(\ref{inl.2}).
\end{proof}

Put $N=N_{1}\cap N_{2}$, where
\[
N_{1}=\{(s,x)\in Q_{\hat{T}}; P_{s,x}(\int_{s}^{T}
|\bar{f}(t,X_t)|^2\,dt<\infty)=1\}^{c},
\]
\[
N_{2}=\{(s,x)\in Q_{\hat{T}};\, \mbox{p.v.\,-}
\int_{s}^{T}\bar{f}(t,X_t)\,d\alpha^{s,x}_{s,t}\,\, \mbox{ exists
and is finite}\,\, P_{s,x}\mbox{-a.s.}\}^{c}
\]
and
\[
\mbox{p.v.\,-}\int_{s}^{T} \bar{f}(t,X_t)\,d\alpha^{s,x}_{s,t}
\equiv\lim_{\delta\rightarrow 0^{+}}
\int_{s+\delta}^{T}\bar{f}(t,X_t)\,d\alpha^{s,x}_{s,t}.
\]

\begin{corollary}\label{coll.2}
If $f\in\mathbb{L}_{2,\varrho}(Q_{\hat{T}})$ then \mbox{\rm
cap}$(N)=0$.
\end{corollary}
\begin{proof}
Follows directly from Proposition \ref{prop.cd} and Proposition
\ref{lml.1}.
\end{proof}

Let $\bar{f}\in {\mathbb{L}}_{2,\varrho}(Q_{T})$. For fixed
$(s,x)\in Q_{\hat{T}}$ we set
\begin{equation}\label{defl.2}
\int_{r}^{t}\bar{f}(\theta,X_{\theta})\,d^{*}X_{\theta}
=-\int_{r}^{t}\bar{f}
(\theta,X_{\theta})(dM_{s,\theta}+d\alpha^{s,x}_{s,\theta})
-\int^{T+s-r}_{T+s-t}
\bar{f}(\bar{\theta},\bar{X}_{\theta})\,dN^{s,x}_{s,\theta}
\end{equation}
for all $0\le s <r \le t\le T$, and for fixed $(s,x)\in N^{c}$ we
set
\begin{align}\label{defl.1}
\int_{s}^{t}\bar{f}(\theta,X_{\theta})\,d^{*}X_{\theta}&=
-\int_{s}^{t}\bar{f}(\theta,X_{\theta})\,dM_{s,\theta}\nonumber\\
&\quad+\mbox{p.v.\,-}\int_{s}^{t}
\bar{f}(\theta,X_{\theta})\,d\alpha^{s,x}_{s,\theta}
-\int^{T}_{T+s-t}
\bar{f}(\bar{\theta},\bar{X}_{\theta})dN^{s,x}_{s,\theta}
\end{align}
for all $0\le s\le t\le T$.

Under stronger integrability conditions on $f$ all integrals on
the right-hand side of (\ref{eq2.6}) are defined for q.e $(s,x)$.

\begin{proposition}
\label{propl.7} If $f\in{\mathbb{L}}_{p,\varrho}(Q_{T})$ for some
$p>2$ then for q.e. $(s,x)\in Q_{\hat{T}}$\,,
\begin{align*}
\int_{s}^{t}\bar{f}(\theta,X_{\theta})\,d^{*}X_{\theta}&=
-\int_{s}^{t}\bar{f}
(\theta,X_{\theta})(dM_{s,\theta}+d\alpha^{s,x}_{s,\theta})\\
&\quad -\int^{T}_{T+s-t}
\bar{f}(\bar{\theta},\bar{X}_{\theta})\,dN^{s,x}_{s,\theta},\quad
t\in[s,T],\quad P_{s,x}\mbox{-}a.s..
\end{align*}
\end{proposition}
\begin{proof}
By (\ref{nq7.a}),
\[
E_{s,x}\int_{s}^{T}|\bar{f}(t,X_t)|\,d|\alpha^{s,x}_{s,\cdot}|_t
\le C(p)(E_{s,x}\int_{s}^{T}
|\bar{f}(t,X_t)|^{p}\,dt)^{2/p}.
\]
From Proposition \ref{prop3.2} it follows that the right-hand side
is finite for a.e. $(s,x)\in Q_{\hat T}$. Hence, by Corollary
\ref{coll.1}, it is finite for q.e.  $(s,x)\in Q_{\hat T}$. The
result now follows from Corollary \ref{coll.2}.
\end{proof}

\begin{proposition}
Let $\{\bar{f}_{n}\}\subset \mathbb{L}_{2,\varrho}(Q_{T})$ and
$\bar{f}_{n}\rightarrow \bar{f}$ in
$\mathbb{L}_{2,\varrho}(Q_{T})$. Then
\begin{enumerate}
\item [\rm(i)] For every $(s,x)\in Q_{\hat{T}}$ and $r\in(s,T]$,
\[E_{s,x}\sup_{r \le t \le T}|
\int_{r}^{t}\bar{f}_{n}(\theta,X_{\theta})\,d^{*}X_{\theta}
-\int_{r}^{t}\bar{f}(\theta,X_{\theta})\,d^{*}X_{\theta}|\rightarrow 0,
\]
\item [\rm(ii)] There exists subsequence (still denoted by $\{n\}$)
such that for q.e. $(s,x)\in Q_{\hat{T}}$\,,
\[
E_{s,x}\sup_{s \le t \le T}|
\int_{s}^{t}\bar{f}_{n}(\theta,X_{\theta})\,d^{*}X_{\theta}
-\int_{s}^{t}\bar{f}(\theta,X_{\theta})\,d^{*}X_{\theta}|\rightarrow 0,
\]
\end{enumerate}
\end{proposition}
\begin{proof}
(i) follows easily from (\ref{inl.1}) and (\ref{inl.2}). (ii)
follows from Proposition \ref{propl.1}, because
\begin{equation*}
\int_{Q_{T}}\left(E_{s,x}\sup_{s\le t \le T}|
\int_{s}^{t}\bar{f}_{n}(\theta,X_{\theta})\,d^{*}X_{\theta}
-\int_{s}^{t}\bar{f}(\theta,X_{\theta})\,d^{*}X_{\theta}|\right)\varrho(x)\,dx
\le \|\bar{f}_{n}-\bar{f}\|^{2}_{2,\varrho,T}
\end{equation*}
by Proposition \ref{prop3.2}.
\end{proof}

\nsubsection{Time-inhomogenous additive functionals and Dirichlet
processes} \label{sec4}

In this section we will be concerned with conditions on $u$ under
which the functional $X^u=\{X^u_{s,t}\equiv u(t,X_t)-u(s,X_s);0\le
s\le t\le T\}$ is a Dirichlet process.

Let $\varrho\in\mathcal{R}_{I}$. For $s\in [0,T)$ we set
$P_{s,\varrho}(\cdot)=\int_{\mathbb{R}^{d}}
P_{s,x}(\cdot)\varrho^{2}(x)\,dx$.

\begin{definition}
{\rm  We say that CAF $A$ of finite variation is locally finite
(square-inte\-gra\-ble) if for every $\eta\in C^{+}_{0}(Q_{T})$,
\[
\int_{0}^{T}E_{s,\varrho}\int_{0}^{T}
\eta(t,X_{t})\,d|A_{s,\cdot}|_{t}\,ds<\infty,
\quad \left(\int_{0}^{T}E_{s,\varrho}|A_{s,\cdot}|^{2}_{T}\,ds<\infty\right).
\]
}
\end{definition}

\begin{definition}
{\rm (i) Let $(s,x)\in Q_{\hat{T}}$ and $r\in[s,T]$. We say that a
$\{\mathcal{G}^{s}_{t}\}$-adapted process $Y$ is a continuous
Dirichlet process on $[r,T]$ under $P_{s,x}$ if
\begin{eqnarray}\label{eqd.17}
Y_{t}=M_{t}+A_{t},\quad t\in[r,T],\quad  P_{s,x}\mbox{-a.s.,}
\end{eqnarray}
where $M$ is a continuous
$(\{\mathcal{G}^{s}_t\},P_{s,x})$-square-integrable martingale on
$[r,T]$ and $A$ is a continuous $\{\mathcal{G}^{s}_t\}$-adapted
process on $[r,T]$ such that $\langle A\rangle^T_r=0$ in the sense
of (\ref{eqp.27}). We say that $Y$ is a continuous Dirichlet
process on $(s,T]$ under $P_{s,x}$ if it is continuous Dirichlet
process on $[r,T]$ under $P_{s,x}$ for every $r\in (s,T]$.\\
(ii) Let $\{{\mit\Pi}_m=\{t_0,t_1,\dots,t_{i(m)}\}\}$ be a
sequence of partitions of $[s,T]$ whose mesh-size converges to
zero as $m\rightarrow\infty$. If $Y$ admits decomposition of the
form (\ref{eqd.17}) with a continuous
$\{\mathcal{G}^{s}_t\}$-adapted process $A$ on $[r,T]$ such that
$\langle A\rangle^T_r=0$ along $\{{\mit\Pi}_m\}$ then we call it a
continuous Dirichlet process along $\{\mit\Pi_{m}\}$. }
\end{definition}

Given $u\in{\mathcal{W}}_{\varrho}$ set
\begin{equation}
\label{eq3.25} M^u_{s,t}\equiv \int^t_s\nabla
u(\theta,X_{\theta})\,dM_{s,\theta}=\int^t_s\sigma\nabla
u(\theta,X_{\theta})\,dB_{s,\theta},\quad 0\le s\le t\le T,
\end{equation}
where $B$ is defined by (\ref{eq2.5}).

\begin{theorem}\label{tw3.1}
Assume that $u\in {\mathcal{W}}_{\varrho}$. Then there exists a
quasi-continuous version of $u$ (still denoted by $u$) such that
\begin{enumerate}
\item[\rm(i)]
For every $(s,x)\in Q_{\hat{T}}$ the functional $X^{u}$ is a
continuous Dirichlet process on $(s,T]$ under $P_{s,x}$ with the
decomposition
\begin{equation}
\label{eq6.4} X^u_{r,t}=M^u_{r,t}+A^u_{r,t},\quad s<r\le t\le
T,\quad P_{s,x}\mbox{-}a.s.,
\end{equation}
where
\begin{equation}
\label{eq6.5}
A^u_{r,t}=\int_{r}^{t}f^{0}(\theta,X_{\theta})\,d\theta +
\int_{r}^{t}a^{-1}\bar{f}(\theta,X_{\theta})\,d^{*}X_{\theta},\quad
s<r\le t\le T,\quad P_{s,x}\mbox{-}a.s.
\end{equation}
with $f^{0},\bar{f}\in\mathbb{L}_{2,\varrho}(Q_{T})$ such that
${\mathcal{L}}u=f^{0}+\dyw\bar{f}$, where
${\mathcal{L}}=\frac{\partial}{\partial t}+L_t$.
\item[\rm(ii)] For q.e. $(s,x)\in Q_{\hat{T}}$ decomposition (\ref{eq6.4}),
(\ref{eq6.5}) holds true with $r=s$.
\item[\rm(iii)] If $\frac{\partial u}{\partial t}
\in \mathbb{L}_{2,\varrho}(Q_{T})
+\mathbb{L}_{2}(0,T;W^{-1}_{p,\varrho})$ with $p>2$ then for q.e.
$(s,x)\in Q_{\hat{T}}$, $X^{u}$ is a Dirichlet process on $[s,T]$
under $P_{s,x}$ with decomposition (\ref{eq6.4}) for $r=s$.
\item[\rm(iv)]For every sequence $\{\mit\Pi_{m}\}$ of partitions of $[0,T]$
whose mesh-size converges to zero as $m\rightarrow\infty$ there
exists a subsequence $\{\mit\Pi_{m'}\}$ such that $X^{u}$ is a
continuous Dirichlet process on $[s,T]$ along
$\{\mit\Pi^{s,T}_{m'}=\mit\Pi_m\cap[s,T]\}$ for q.e. $(s,x)\in
Q_{\hat{T}}$ admitting decomposition (\ref{eq6.4}) for $r=s$.
\end{enumerate}
\end{theorem}
\begin{proof}
(i) First note that it is known that if $u\in
\mathcal{W}_{\varrho}$ then there exist
$f^{0},\bar{f}\in\mathbb{L}_{2,\varrho}(Q_{T})$ such that
${\mathcal{L}}u=f^{0}+\dyw\bar{f}$. Let us fix $(s,x)\in Q_{\hat
T}$. Let
$\Phi^\varepsilon=f^0_\varepsilon+\dyw{\bar{f}_{\varepsilon}}$,
where $f^i_\varepsilon$, $i=0,\dots,d$, are standard
mollifications of $f^i$, and let
$\bar{f}_{\varepsilon}=(f^{1}_{\varepsilon},\dots,f^{d}_{\varepsilon})$.
Let $u^{\varepsilon}(T)$ be the standard mollification of $u(T)$ and  $u_n$ be a continuous version on $Q_{T}$ of weak solution of
the Cauchy problem
\begin{equation}
\label{eq3.14} (\frac{\partial}{\partial t}+L_t)u_{n}=\Phi^n,\quad
u_{n}(T)=u^{n}(T),
\end{equation}
where $\Phi^n=\Phi^{\varepsilon},\,\, u^{n}(T)=u^{\varepsilon}$
with $\varepsilon=1/n$. From \cite{A.Roz.BSDE} we know that
$P_{s,x}$-a.s.,
\begin{equation}\label{eq}
X^{u_n}_{r,t}=M^{u_n}_{r,t}+A^{u_n}_{r,t},
\end{equation}
for $s\le r \le t \le T$, where
\[
A^{u_n}_{r,t}=\int_{r}^{t}(f_n^{0}
+\mbox{div}\bar{f}_n)(\theta,X_\theta)\,d\theta
\]
Let us define a version of $u$ (still denoted by $u$) as follows:
$u(s,x)=\lim_{n\rightarrow\infty} u_{n}(s,x)$ if the limit exists
and zero otherwise. It is known (see \cite{Lad}) that $u_{n}\rightarrow u$ in $\mathcal{W}_{\varrho}$.
Next, let
\[
A^{u}_{r,t}=\int_{r}^{t}f^{0}(\theta,X_{\theta})\,d\theta
+\int_{r}^{t}a^{-1}\bar{f}(\theta,X_\theta)\,d^{*}X_{\theta},\quad
s\le r \le t \le T.
\]
By Proposition \ref{propl.1} and (\ref{inl.1}),  for every
$(s,x)\in Q_{\hat{T}}$ and $r\in (s,T]$, $X^{u_n}_{r,t}\rightarrow
X^{u}_{r,t}$, $M^{u_n}_{r,t}\rightarrow M^{u}_{r,t}$,
$A^{u_n}_{r,t}\rightarrow A^{u}_{r,t}$ in $\mathbb{L}_{1}(\Omega,
P_{s,x})$ uniformly in $t\in[r,T]$. Therefore passing to the limit
in (\ref{eq}) we get (\ref{eq6.4}), (\ref{eq6.5}). By
(\ref{defl.2}),  for every $(s,x)\in Q_{\hat{T}}$ and $r\in
(s,T]$,
\[
\left\langle
\int_{r}^{\cdot}
a^{-1}\bar{f}(\theta,X_{\theta})\,d^{*}X_{\theta}\right\rangle_{r}^{T}
=\left\langle
\int_{r}^{\cdot}a^{-1}\bar{f}(\theta,X_{\theta})\,dM_{s,\theta}
+\int^{T+s-r}_{T+s-\cdot}a^{-1}\bar{f}
(\bar{\theta},\bar{X}_{\theta})\,dN^{s,x}_{s,\theta}\right\rangle_{r}^{T}
\]
under  $P_{s,x}$. Let $\{\bar{f}_{n}\}\subset
C_{0}^{\infty}(Q_{T})$ be a sequence such that
$\bar{f}_{n}\rightarrow \bar{f}$ in
${\mathbb{L}}_{2,\varrho}(Q_{T})$. Then
\begin{align}\label{sts}
&E_{s,x}\sum_{t_i\in\mit\Pi_m^{r,T}}\left|\int_{t_{1}}^{t_{2}}
(\bar{f}-\bar{f}_{n})(\theta,X_{\theta})\,dM_{s,\theta}
+\int^{T+s-t_{1}}_{T+s-t_{2}}(\bar{f}-\bar{f}_{n})
(\bar{\theta},\bar{X}_{\theta})\,dN^{s,x}_{s,\theta}\right|^{2}\nonumber\\
&\qquad\le CE_{s,x}\int_{r}^{T}|\bar{f}_{n}-\bar{f}|^{2}
(\theta,X_{\theta})\,d\theta \le C\varrho^{-2}(x)(r-s )^{-d/2}
\|\bar{f}_{n}-\bar{f}\|_{2,\varrho,T}\,.
\end{align}
From this and the fact that $\langle
\int_{r}^{\cdot}\bar{f_n}(\theta,X_{\theta})\,
d^{*}X_{\theta}\rangle_{r}^T=0,\, r\in(s,T],\, P_{s,x}$-a.s. for
every $(s,x)\in Q_{\hat{T}}$ we get the first assertion of (i). To
prove  (ii) it suffices to pass to the limit with $r\rightarrow
s^{+}$ in (\ref{eq6.4}) (if the limit exists) and use Corollary
\ref{coll.1}. Since from (i) it follows that $\langle
\int_{r}^{\cdot}\bar{f}(\theta,X_{\theta})\,
d^{*}X_{\theta}\rangle_{r}^T=0,\, r\in(s,T],\, P_{s,x}$-a.s. for
every $(s,x)\in Q_{\hat{T}}$, to prove (iii) it suffices to show
that for q.e $(s,x)\in Q_{\hat{T}}$ there exists the covariation
$\left\langle \int_{s}^{\cdot}
a^{-1}\bar{f}(\theta,X_{\theta})\,d^{*}X_{\theta}\right\rangle_{s}^{T}$
under $P_{s,x}$. But the last statement is a direct consequence of
Proposition \ref{propl.7}. Finally, to prove (iv) let us set
$B=\{(s,x)\in Q_{T}:\limsup_{m\rightarrow\infty}
E_{s,x}\sum_{t_{i}\in\mit\Pi_{m}}|A_{s,t_{i+1}}-A_{s,t_{i}}|^{2}>0\}$
and $\tau=\inf\{t\in [s,T]:(t,X_{t})\in K\}$, where $K$ is a
compact subset such that $K\subset B$. Then by strong Markov
property with random shift and additivity of $A$,
\begin{align*}
P_{s,x}(\tau<\infty)& =P_{s,x}(\limsup_{m\rightarrow\infty}
E_{\tau\wedge T,X_{\tau\wedge T}} \sum_{t_{i}\in\mit\Pi_{m}}
|A_{\tau\wedge T,t_{i+1}}-A_{\tau\wedge T,t_{i}}|^{2}>0)\\
&= P_{s,x}(\limsup_{m\rightarrow\infty}E_{s,x}(\sum_{t_{i}\in
\mit\Pi_{m}}|A_{\tau\wedge T,t_{i+1}}-A_{\tau\wedge T,t_{i}}|^{2})
|\mathcal{G}^{s}_{\tau\wedge T})>0)\\
&\le P_{s,x}(\limsup_{m\rightarrow\infty}E_{s,x}(\sum_{t_{i}\in
\mit\Pi_{n}}|A_{s,t_{i+1}}-A_{s,t_{i}}|^{2})
|\mathcal{G}^{s}_{\tau\wedge T})>0).
\end{align*}
Set
\[
T_{m}(s,x,\omega) =E_{s,x}(\sum_{t_{i}\in
\mit\Pi_{m}}|A_{s,t_{i+1}}-A_{s,t_{i}}|^{2})
|\mathcal{G}^{s}_{\tau\wedge T})
\]
and define the measure $\Pi$ as in the proof of  Proposition
\ref{prop.cd}. Since we know already that $\langle
A_{s,\cdot}\rangle^T_s=0$ under $P_{s,x}$ for a.e. $(s,x)\in Q_T$,
it follows that $T_{n}\wedge M\rightarrow 0$ in
${\mathbb{L}_1}(Q_{T}\times\Omega,\Pi)$  and hence that there
exists a subsequence (still denoted by $m$) such that
$T_{m}\rightarrow 0$, $\Pi$-a.e.. Therefore $T_{m}(s,x)\rightarrow
0$, $P_{s,x}$-a.s. for a.e. $(s,x)\in Q_{T}$, which proves that
$\mbox{cap}_{L}(K)=0$, hence that $\mbox{cap}_{L}(B)=0$ by Remark
\ref{choquet}.
\end{proof}

\begin{corollary}
For every $\Phi\in \mathbb{L}_{2}(0,T;H^{-1}_{\varrho})$ there
exists a unique CAF $A$ of zero quadratic variation such that
\[
A_{r,t}=\int_{r}^{t}f^{0}(\theta,X_{\theta})\,d\theta +
\int_{r}^{t}a^{-1}\bar{f}(\theta,X_{\theta})\,d^{*}X_{\theta},\quad
s<r\le t\le T,\quad P_{s,x}\mbox{-}a.s.
\]
for any decomposition of $\Phi$ of the form
$\Phi=f^{0}+\dyw\bar{f}$,
$f^{0},\bar{f}\in\mathbb{L}_{2,\varrho}(Q_{T})$.

\end{corollary}

In the sequel we write
$\int_{r}^{t}d\Phi(\theta,X_{\theta})=A_{r,t}$, $s\le r\le t\le T$
or $\Phi\sim A$ if $A$ is the CAF corresponding to
$\Phi\in\mathbb{L}_{2}(0,T;H^{-1}_{\varrho})$ in the sense of the
above corollary.

\begin{remark}\label{linear}
{\rm From linearity of the operator $L_{t}$ it follows immediately
that the mapping $\mathcal{W}_{\varrho}\ni u\rightarrow A^{u}$,
where $A^u$ is the functional of Theorem \ref{tw3.1}, is linear. }
\end{remark}

It is worth mentioning that the decomposition (\ref{eq6.4})
implies Fukushima's decomposition of $X^u$ into martingale AF of
finite energy and CAF of zero energy (for related results for
time-independent $u$ see \cite{Roz.dec}). To state the result, let
us recall first the definition of energy of time-inhomogeneous
additive functionals of $\BX$ and its basic properties.

\begin{definition}
{\rm Let $A,B$ be CAFs of $\BX$. We define the mutual energy of
$A$ and $B$ by
\[
e(A,B)=\lim_{h\rightarrow 0^{+}}\int_{0}^{T-h}
E_{s,\varrho}A_{s,s+h}B_{s,s+h}\,ds
\]
(whenever the limit exists), and we put  $e(A)=e(A,A)$. }
\end{definition}

One can check that the energy has the following properties.
\begin{enumerate}
\item [(i)] $e(A+B)=e(A)+e(B)+2e(A,B)\le 2(e(A)+e(B))$,
\item [(ii)] If $e(A)=0$ and $e(B)<\infty$, then $e(A,B)=0$ and
$e(A+B)=e(B)$.
\end{enumerate}

\begin{lemma}\label{lm.76}
If $\Phi\in \mathbb{L}_{2}(0,T;H^{-1}_{\varrho})$ then
\[
e(\int d\Phi(\cdot,X_{\cdot}))\le C\|\Phi\|^{2}_{*}\,.
\]
\end{lemma}
\begin{proof}
Let $s\in (0,T)$,  $h\in(0,T-s)$ and let
$u\in\mathcal{W}_{\varrho}$ be a solution of PDE$(0,\Phi)$ on
$[s,s+h]$.  By Theorem \ref{tw3.1} and Remark \ref{rem.76} there
exists a quasi-continuous version of $u$ (still denoted $u$) such
that for a.e.  $x\in \mathbb{R}^{d}$,
\[
X^u_{s,t}=M^u_{s,t}+A^u_{s,t},\quad t\in[s,s+h],\quad
P_{s,x}\mbox{-}a.s.,
\]
where
\begin{equation*}
A^u_{s,t}= \int_{s}^{t}d\Phi(\theta,X_{\theta}),\quad
t\in[s,s+h],\quad P_{s,x}\mbox{-}a.s..
\end{equation*}
Hence, by Proposition \ref{prop3.2} and Theorem \ref{tha.1},
\begin{align*}
E_{s,\varrho}|\int_{s}^{s+h} d\Phi(\theta,X_{\theta})|^2 &\le
C(\|\nabla u\|^{2}_{2,\varrho,s,s+h}+\sup_{s\le t\le s+h}
\|u(t)\|^{2}_{2,\varrho})\\& \le
C(\|f^{0}\|_{2,\varrho,s,s+h}+\|\bar{f}\|_{2,\varrho,s,s+h}),
\end{align*}
where $f^{0},\bar{f}\in\mathbb{L}_{2,\varrho}(Q_{T})$ are such
that $\Phi=f^{0}+\dyw(\bar{f})$. From the above inequality  the
result easily follows.
\end{proof}

\begin{corollary}
Let $M^u,A^u$ be AFs of the decomposition (\ref{eq6.4}). Then
$e(M^{u})<\infty$, $e(A^{u})=0$.
\end{corollary}
\begin{proof}
Using (\ref{eq3.25}) and Proposition \ref{prop3.2} one can check
that $e(M^u)\le C\|\nabla u\|_{2,\varrho,T}^{2}<\infty$. To prove
that $e(A^u)=0$ let us write $\mathcal{L}u=\Phi$ and define
$\Phi^{n}$, $A^{u_{n}}$ as in the proof of Theorem \ref{tw3.1}.
Since the CAF $A^{u_{n}}$ has finite variation, direct calculation
shows that $e(A^{u_{n}})=0$. From this, Lemma \ref{lm.76}
and property (i), it follows that
\[
e(A^u)\le 2e(A^u-A^{u_n})\le
C(\|f^{0}-f^{0}_{n}\|^{2}_{2,\varrho,T}
+\|\bar{f}-\bar{f}_{n}\|^{2}_{2,\varrho,T})
\]
for $n\in\mathbb{N}$ which completes the proof.
\end{proof}

\nsubsection{Continuous additive functionals of zero-quadratic
variation} \label{sec5}

Given CAF $A$ we set
\[
D(A)=\{(s,x)\in Q_{\hat{T}}: A_{s,\cdot}\mbox{
is continuous and additive on }[s,T]\mbox{ under }P_{s,x}\}.
\]
and
\[
D_{0}(A)=\{(s,x)\in D(A); E_{s,x}\int_{s}^{T}|A_{s,t}|^{2}\,dt
<\infty\}.
\]

For instance, if
$A_{s,t}=\int_{s}^{t}f(\theta,X_{\theta})\,d\theta$, $0\le s\le
t\le T$, for some $f\in\mathbb{L}_{1,\varrho}(Q_{T})$ then
$\{(s,x)\in Q_{\hat{T}};
P_{s,x}(\int_{s}^{T}|f(t,X_t)|\,dt<\infty)=1\}\subset D(A)$, and if
\begin{equation}
\label{eq5.1}
A_{s,t}=\int_s^t\bar{f}(\theta,X_{\theta})\,d^{*}X_{\theta}, \quad
0\le s\le t\le T,
\end{equation}
for some $\bar{f}\in\mathbb{L}_{2,\varrho}(Q_{T})$ then
$N^{c}\subset D(A)$, where $N$ is defined in Corollary \ref{coll.2}.

\begin{lemma}
Let $\bar{f}\in\mathbb{L}_{2,\varrho}(Q_{T})$ and let $A$ be
defined by (\ref{eq5.1}). Then \mbox{\rm
cap}$_{L}((D_{0}(A))^{c})=0$.
\end{lemma}
\begin{proof}
Let $u\in\mathcal{W}_{\varrho}$ be a solution of PDE$(0,\Phi)$,
where $\Phi=\dyw(a^{-1}\bar{f})$ By Theorem \ref{tw3.1} and Remark
\ref{rem.76} there exists a quasi-continuous version of $u$ (still
denoted $u$) such that for for every $s\in (0,T)$ and a.e.  $x\in
\mathbb{R}^{d}$,
\[
X^u_{s,t}=M^u_{s,t}+A^u_{s,t},\quad s\le t\le T,\quad
P_{s,x}\mbox{-}a.s.,
\]
where
\[
A^u_{s,t}= \int_{s}^{t}\bar{f}(\theta,X_{\theta})\,d^{*}X_{\theta}\quad
s\le t\le T,\quad P_{s,x}\mbox{-}a.s..
\]
Hence, by Proposition \ref{prop3.2} and Theorem \ref{tha.1},
\[
E_{s,\varrho}A_{s,t}^{2}\le C\|\Phi\|_{*}, \quad 0<s\le t\le T.
\]
Consequently, $E_{s,x}\int_{s}^{T}|A_{s,t}|^{2}\,dt<\infty$ for
a.e. $(s,x)\in Q_{\hat T}$ and hence for  q.e. $(s,x)\in Q_{\hat
T}$ by Corollary \ref{coll.1}.
\end{proof}

Let $\Phi\in {\mathbb{L}}_{2}(0,T;H^{-1}_{\varrho})$, $A\sim\Phi$
and let $f^{0},\bar{f}\in\mathbb{L}_{2,\varrho}(Q_{T})$ be such
that $\Phi=f^{0}+\dyw\bar{f}$. Our next goal is to define the
integral with respect to $A$ and show that $A$ is determined by
its $\alpha$-potential.

Given $\eta\in{\mathcal W}_{\varrho}\cap \mathcal{B}_b(Q_{T})$ we
set
\begin{align}\label{di}
(\eta\cdot A)_{r,t}=
\int_{r}^{t}\eta(\theta,X_{\theta})\,dA_{s,\theta}
&\equiv\int_{r}^{t}\eta f^{0}(\theta,X_{\theta})\,d\theta
-\int_{r}^{t}\nabla\eta\bar{f}(\theta,X_{\theta})\,d\theta\nonumber\\
&\quad+\int_{r}^{t}a^{-1}\bar{f}\eta(\theta,X_{\theta})\,d^{*}X_{\theta},
\quad 0\le s<r \le t\le T.
\end{align}
Observe that  from (\ref{inl.1}), (\ref{inl.2}) it follows that
all the integrals on the right-hand side of (\ref{di})  are well
defined. Moreover, setting
\[
N_{2}= (D(\int \eta f^{0}\,dt))^{c}\cup
(D(\int\nabla\eta\bar{f}\,dt))^{c}\cup (D(\int
a^{-1}\bar{f}\eta\,d^{*}X))^{c}
\]
we see that cap$_{L}(N_{2})=0$ and for every $(s,x)\in N_{2}^{c}$
the right-hand side of (\ref{di}) converges $P_{s,x}$-a.s. to a
finite limit as  $r\rightarrow s^{+}$. Thus, (\ref{di}) defines a
CAF of $\BX$.

From the following proposition it follows in particular that
$\eta\cdot A$ does not depend on the choice of $f^{0},\bar{f}$ in
the decomposition of $\Phi$.

\begin{proposition}\label{apr.1}
Let $\Phi\ni\mathbb{L}_{2}(0,T,H^{-1}_{\varrho})$ and let
$A\sim\Phi$.
\begin{enumerate}
\item[\rm(i)]For every bounded
$\eta\in{\mathcal{W}}_{\varrho}$ there exists a sequence
$\{A^{n}\}$ of locally finite CAFs of finite variation such that
for q.e. $(s,x)\in Q_{\hat T}$,
\[
E_{s,x}\sup_{s\le t \le
T}|\int_{s}^{t}\eta(\theta,X_{\theta})\,dA^{n}_{s,\theta}
-\int_{s}^{t}\eta(\theta,X_{\theta})\,dA_{s,\theta}|\rightarrow 0
\]
\item[\rm(ii)]There exists a sequence $\{A^{n}\}$ of locally finite
CAFs of finite variation such that for every bounded
$\eta\in{\mathcal{W}}_{\varrho}$, $(s,x)\in Q_{\hat{T}}$ and
$r\in(s,T]$,
\[
E_{s,x}\sup_{r\le t \le T}
|\int_{r}^{t}\eta(\theta,X_{\theta})\,dA^{n}_{r,\theta}
-\int_{r}^{t}\eta(\theta,X_{\theta})\,dA_{r,\theta}|\rightarrow
0.
\]
\end{enumerate}
\end{proposition}
\begin{proof}
Let $A^n=A^{u_n}$, where $A^{u_n}$ is defined as in the proof of
Theorem \ref{tw3.1}. Then the second part follows immediately from
the definition of $\eta\cdot A^{n}$, $\eta\cdot A$ and
(\ref{inl.1}), (\ref{inl.2}). To prove the first part, let us
observe that by Proposition \ref{prop3.2},
\[
\int_{Q_{T}}(E_{s,x}\sup_{s\le t \le T}|(\eta\cdot A^{n})_{s,t}
-(\eta\cdot A)_{s,t}|)\varrho(x)\,dx \rightarrow 0,
\]
so the result follows from Proposition \ref{propl.1}.
\end{proof}

\begin{remark}
{\rm Notice that from Proposition \ref{apr.1} it follows that if
$A$ is a CAF of finite variation corresponding to some $\Phi\in
{\mathbb{L}}_{2}(0,T;H^{-1}_{\varrho})$ then the usual
Lebesgue-Stieltjes integral $\int_s^{\cdot}\eta(t,X_t)\,dA_{s,t}$
and the integral in the sense of (\ref{di}) coincide. }
\end{remark}

Using the  definition (\ref{di}) of the integral
with respect to additive functionals of zero-quadratic variation
we can define Laplace transform of such an additive functional.

For $\alpha>0$ we put
\[
U^{\alpha}_{A}(s,x) =E_{s,x}\int_{s}^{T}e^{-\alpha
(t-s)}\,dA_{s,t},\quad (s,x)\in D(e^{-(\cdot-s)}\cdot A)
\]
and
\begin{equation}
\label{eq4.03} U_{A}^{\alpha}\eta(s,x)
=E_{s,x}\int_{s}^{T}e^{-\alpha(t-s)}\eta(t,X_{t})\,dA_{s,t},\quad
(s,x)\in D(e^{-(\cdot-s)}\eta\cdot A)
\end{equation}
for $\eta\in{\mathcal{W}}_{\varrho}\cap C_{b}(Q_{T})$. In case $A$
is a CAF of finite variation, the integral in (\ref{eq4.03}) is
the usual Lebesgue-Stieltjes integral which is well defined for
all $\eta\in C_{b}(Q_{T})$.

In the sequel we denote $D(e^{-(\cdot-s)}\eta\cdot A)$ by
$D(U_{A}^{\alpha}\eta)$. If $A_{s,t}=t-s$, then we denote
$U^{\alpha}_A\eta$ by $U^{\alpha}\eta$.

Let $\{R_{\alpha};\alpha\ge 0\}$ denote the resolvent of
${\mathcal{L}}$ on $\mathbb{L}_{2,\varrho}(Q_T)$. Notice that if
$\xi\in \mathcal{B}_{b}(Q_{T})$ then $R_{\alpha}\xi\in
\mathcal{W}_{\varrho}\cap\mathcal{B}_{b}(Q_{T})$ and $\nabla
R_{\alpha}\xi\in\mathcal{B}_{b}(Q_{T})$. Indeed, the first
assertion follows immediately from the fact  that $R_{\alpha}\xi$
is a strong solution of the Cauchy problem $(\alpha+\LL)u=-\xi$,
$u(T)=0$ and the representation formula
\[
R_{\alpha}\xi(s,x)=E_{s,x}\int_{s}^{T}
e^{-\alpha(\theta-s)}\xi(\theta,X_{\theta})\,d\theta, \quad
(s,x)\in Q_{\hat T}\,.
\]
The second assertion follows from the formula
\[
\nabla R_{\alpha}\xi(s,x)=\int_{Q_{sT}}e^{-\alpha(\theta-s)}
\xi(\theta,y)\nabla_{x}p(s,x,\theta,y)\,d\theta\,dy,\quad (s,x)\in
Q_{\hat T}
\]
and integrability of $\nabla_{x}p(s,x,\cdot,\cdot)$ proved in
\cite[Theorem 10]{Aronson}.

\begin{proposition}
Let $A$ be CAF associated with some
$\Phi\in\mathbb{L}_{2}(0,T;H^{-1}_{\varrho})$.  Then
\[
D_{0}(A)\subset
\bigcap_{\xi\in\mathcal{B}_b(Q_{T}),\alpha,\beta\ge0}
D_{0}(e^{-\alpha(\cdot-s)}\cdot R_{\beta}(\xi)\cdot A).
\]
\end{proposition}
\begin{proof}
Let $\eta\in R_{\beta}(\mathcal{B}_{b}(Q_{T}))$  and $\xi\in
\mathcal{B}_{b}(Q_{T})$ be such that $\eta=R_{\beta}\xi$. By
Proposition \ref{apr.1}, for every $(s,x)\in Q_{\hat{T}}$ and
$r\in(s,T]$,
\[
E_{s,x}\sup_{r\le t \le T} |\int_{r}^{t}e^{-\alpha(\theta-s)}
\eta(\theta,X_{\theta})\,dA^{n}_{r,\theta}
-\int_{r}^{t}e^{-\alpha(\theta-s)}
\eta(\theta,X_{\theta})\,dA_{r,\theta}|\rightarrow0
\]
for some sequence $\{A^{n}\}$ of CAFs of finite variation.  By
results proved in \cite{A.Roz.BSDE} and  elementary calculations,
for every $(s,x)\in Q_{\hat{T}}$ we have
\begin{align*}
\eta(t,X_{t})&=\int_{t}^{T}e^{-\beta(\theta-s)}
\xi(\theta,X_{\theta})\,d\theta\\
&\quad-\int_{t}^{T} e^{-\beta(\theta-s)}
\sigma\nabla\eta(\theta,X_{\theta})\,dB_{s,\theta},\quad
t\in[s,T],\quad P_{s,x}\mbox{-}a.s..
\end{align*}
Hence applying the integration by parts formula to
$A_{r,\cdot}^{n}(e^{-\alpha(\cdot-s)}\eta(\cdot,X_{\cdot}))$ and letting $n\rightarrow\infty$  we
conclude that for every $(s,x)\in Q_{\hat{T}}$ under the measure
$P_{s,x}$\,,
\begin{align}
\label{eq5.6} &\int_{r}^{t}e^{-\alpha(\theta-s)}
\eta(\theta,X_{\theta})\,dA_{r,\theta}\nonumber\\
&\quad=\int_{r}^{t}e^{-(\alpha+\beta)(\theta-s)}
\xi(\theta,X_{\theta})A_{r,\theta}\,d\theta
-\int_{r}^{t}e^{-(\alpha+\beta)(\theta-s)}
\nabla\eta(\theta,X_{\theta})A_{r,\theta}\,dB_{s,\theta}\nonumber\\
&\qquad-\alpha\int_{r}^{t}
e^{-(\alpha+\beta)(\theta-s)}\eta(\theta,X_{\theta})
A_{r,\theta}\,d\theta +\eta(t,X_{t})e^{-\alpha(t-s)}A_{r,t}
\end{align}
for $s<r\le t\le T$. Since
$\eta,\nabla\eta\in{\mathcal{B}}_b(Q_T)$ and $r\mapsto\eta(r,X_r)$
is continuous, letting $r\rightarrow s^+$ we get (\ref{eq5.6}) for
$(s,x)\in D(A)$. From (\ref{eq5.6}) with $r=s$ the proposition
easily follows.
\end{proof}

\begin{proposition}\label{LT}
Let $A,D$ be CAFs  associated with some functionals in $\mathbb{L}_{2}(0,T,H^{-1}_{\varrho})$
such that $U^{\alpha}_{D}\eta=U^{\alpha}_{A}\eta$ on
$D_{0}(U^{\alpha}_{D}\eta) \cap D_{0}(U^{\alpha}_{A}\eta)$ for
every $\alpha>0$ and $\eta\in {\mathcal{W}}_{\varrho}\cap
C_{c}(Q_{T})$. Then $A=D$.
\end{proposition}
\begin{proof}
Without lost of generality we may assume that
$U^{\alpha}_{D}\eta=U^{\alpha}_{A}\eta$ on $D(U^{\alpha}_{D}\eta)
\cap D(U^{\alpha}_{A}\eta)$  for every
$\eta\in{\mathcal{W}}_{\varrho}\cap C_{b}(Q_{T})$ because we can
consider functionals $f\cdot A, f\cdot D$ with
$f\in{\mathcal{W}}_{\varrho}\cap C_{c}(Q_{T})$ and from the
equality $f\cdot A=f\cdot D$ for every such $f$ one can deduce
that $A=D$. First we show that
\[
U^{\alpha}_{A}(U^{\alpha}\eta)(s,x)
=E_{s,x}\int_{s}^{T}e^{-\alpha(t-s)}\eta(t,X_{t})A_{s,t}\,dt,\quad
(s,x)\in D_{0}(A).
\]
It is well known (see \cite{Lad}) that $U^{\alpha}\eta\in
{\mathcal{W}}_{\varrho}\cap C_{b}(Q_{T})$, so the above equality
makes sense. Using the Markov property, Proposition \ref{apr.1}
and Fubini's theorem we have that for every $(s,x)\in Q_{\hat{T}}$ and $r\in (s,T]$
\begin{align*}
&E_{s,x}\int_{r}^{T}e^{-\alpha(t-s)}U^{\alpha}\eta(t,X_{t})\,dA_{r,t}
=\lim_{n\rightarrow\infty}E_{s,x}\int_{r}^{T}e^{-\alpha(t-s)}
U^{\alpha}\eta(t,X_{t})\,dA^{n}_{r,t}\\
&\quad= \lim_{n\rightarrow\infty}
E_{s,x}\int_{r}^{T}e^{-\alpha(t-s)}
\left(E_{t,X_{t}}\int_{t}^{T}e^{-\alpha(\theta-t)}
\eta(\theta,X_{\theta})\,d\theta\right)dA^{n}_{r,t}\\
&\quad=\lim_{n\rightarrow\infty}E_{s,x}\int_{r}^{T}e^{-\alpha(t-s)}
E_{s,x}\left(\int_{t}^{T}e^{-\alpha(\theta-t)}
\eta(\theta,X_{\theta})\,d\theta|{\mathcal{G}}^{s}_{t}\right)\,dA^{n}_{r,t}
\\&\quad=\lim_{n\rightarrow\infty}E_{s,x}\int_{r}^{T} e^{-\alpha(t-s)}
\eta(t,X_t)A^{n}_{r,t}\,dt=E_{s,x}\int_{r}^{T} e^{-\alpha(t-s)}
\eta(t,X_t)A_{r,t}\,dt.
\end{align*}
Passing to the limit with $r\rightarrow s^{+}$ for every $(s,x)\in
D_{0}(A)$ we get that
\[
U_{A}^{\alpha}(U^{\alpha}\eta)(s,x)=E_{s,x}\int_{s}^{T}
e^{-\alpha(t-s)} \eta(t,X_t)A_{s,t}\,dt,\quad (s,x)\in D_{0}(A).
\]
By the above and the assumptions it follows that
\[
E_{s,x}\int_{s}^{T}e^{-\alpha(t-s)}\eta(t,X_{t})A_{s,t}\,dt =
E_{s,x}\int_{s}^{T}e^{-\alpha(t-s)} \eta(t,X_{t})D_{s,t}\,dt,
\quad (s,x)\in D_{0}(A).
\]
Hence
\[
E_{s,x}\eta(t,X_{t})A_{s,t}= E_{s,x}\eta(t,X_{t})D_{s,t},\quad
t\in[s,T]
\]
for $(s,x)\in D_{0}(A)$ by the well known properties of the
Laplace transform. Consequently, using the  Markov property and
additivity of $A$, $D$ for every $0\le s\le s^{'}\le t^{'}\le t
\le T$ we have
\begin{align*}
E_{s,x}\eta(t',X_{t'})A_{s',t}
&=E_{s,x}\eta(t',X_{t'})A_{s',t'}+E_{s,x}\eta(t',X_{t'})A_{t',t}\\
&=E_{s,x}(E_{s',X_{s'}}(\eta(t',X_{t'})A_{s',t'}))+
E_{s,x}(\eta(t',X_{t'})E_{t^{'},X_{t'}}A_{t',t})\\
&=E_{s,x}(E_{s',X_{s'}}(\eta(t',X_{t'})D_{s',t'}))
+E_{s,x}(\eta(t',X_{t'})E_{t',X_{t'}}(D_{t',t}))\\
&=E_{s,x}\eta(t',X_{t'})D_{s',t}.
\end{align*}
By induction, we get
\[
E_{s,x}\prod_{i=1}^{k}\eta(t_{i},X_{t_{i}})A_{t',t}
=E_{s,x}\prod_{i=1}^{k}\eta(t_{i},X_{t_{i}})D_{t',t}
\]
for $0\le s\le t'\le t_{1}\le\dots\le t_{k}\le t\le T$ from which
the lemma follows.
\end{proof}

\nsubsection{The semimartingale structure of additive functionals}
\label{sec6}

In this section we proceed with the study of the structure of the
functional $X^u$. We will be concerned with additional conditions
on $u\in {\mathcal{W}}_{\varrho}$ under which $X^u$ is a
semimartingale.

Let $S^{c}$ denote the set of all positive measures on $Q_{T}$
such that $\mu_{|\check{Q}_{T}}\ll\overline{\mbox{cap}}$ and
$\mu(\{0\}\times \mathbb{R}^{d})=\mu(\{T\}\times
\mathbb{R}^{d})=0$, and let $S^{c}_{0}$ be the set of  measures
$\mu\in S^{c}$ for which there exist
$\Phi\in\mathbb{L}_{2}(0,T;H^{-1}_{\varrho})$ such that (\ref{RM})
holds for every $\eta\in C^{\infty}_{c}(Q_{T})$. First we assume
that $\mathcal{L}u\equiv\frac{\partial u}{\partial t}+L_tu\in
S^{c}_{0}-S^{c}_{0}$ and then we consider the case where
$\mathcal{L}u\in \mathcal{M}$. Of course, the first assumption
implies the second one, but in general the converse implication is
not true (see, e.g., \cite[Example I.1]{HJ}). Let us remark also
that in general the functional $\mathcal{L}u$ is not a measure.
For instance, if $d=1$,  and $\mathcal{L}u=f'$, then
$\mathcal{L}u$ is a measure iff $f$ is locally of finite variation
(see, e.g., \cite[Proposition 3.6]{AFP}). Finally, it is worth
noting that the first assumption on the decomposition of
$\mathcal{L}u$ appears naturally when considering obstacle
problems (see, e.g., \cite{MM} and references therein).

\begin{proposition}\label{tw3.2}
Assume that $u\in {\mathcal{W}}_{\varrho}$,  $\mathcal{L}u\in
S^{c}_{0}-S^{c}_{0}$. Then there exist a quasi-continuous version
of $u$ (still denoted by $u$)  and square-integrable positive CAFs
$C, R$ such that for every $(s,x)\in Q_{\hat{T}}$,
\begin{equation}
\label{eq6.1} X^u_{r,t}=M^{u}_{r,t}+C_{r,t}-R_{r,t}, \quad 0
\le s<r\le t\le T\quad P_{s,x}\mbox{-}a.s.,
\end{equation}
\begin{equation}
\label{eq6.1.5} E_{s,x}|C_{r,T}|^2 \le
C\frac{\varrho^{2}(x)}{(r-s)^{d/2}}\|\mu_{1}\|_{\ast}\,, \quad
E_{s,x}|R_{r,T}|^2 \le
C\frac{\varrho^{2}(x)}{(r-s)^{d/2}}\|\mu_{2}\|_{*}
\end{equation}
and
\begin{align}
\label{eq6.2}E_{s,x}\int_r^T\xi(t,X_t)\,dC_{r,t}
&=\int_{Q_{rT}}\xi(t,y)p(s,x,t,y)\,d\mu_1(t,y),\\
\label{eq6.3} E_{s,x}\int_r^T\xi(t,X_t)\,dR_{r,t}
&=\int_{Q_{rT}}\xi(t,y)p(s,x,t,y)\,d\mu_2(t,y)
\end{align}
for all $\xi\in C_{0}(Q_{T})$, where $\mu_1, \mu_2\in S^{c}_{0}$
are such that  $\mathcal{L}u=\mu_{1}-\mu_{2}$. Moreover, for q.e.
$(s,x)\in Q_{\hat{T}}$, (\ref{eq6.1}), (\ref{eq6.2}),
(\ref{eq6.3}) hold with $r=s$.
\end{proposition}
\begin{proof}
By Theorem \ref{tw3.1} there exists CAF $A^{u}$ such that
(\ref{eq6.4}) holds. We are going to show that $A^{u}$ is a CAF of
finite variation. Since $\mu_{1},\mu_{2}\in S^{c}_{0}$, there
exist
$f^{0},g^{0},\bar{f},\bar{g}\in{\mathbb{L}}_{2,\varrho}(Q_{T})$
such that $\mu_{1}=f^{0}+\dyw{\bar{f}}$,
$\mu_{2}=g^{0}+\dyw{\bar{g}}$. Let
$f^{0}_{n},g^{0}_{n},\bar{f}_{n},\bar{g}_{n}$ denote standard
mollifications of  $f^{0},g^{0},\bar{f},\bar{g}$, respectively,
and  let $\mu_{1}^{n}=f^{0}_{n}+\dyw{\bar{f}_{n}}$,
$\mu_{2}^{n}=g_{n}^{0}+\dyw{\bar{g}_{n}}$. It is clear that
$\mu_{1}^{n},\mu_{2}^{n}$ are positive and
$\mu_{1}^{n},\mu_{2}^{n}\in \mathbb{L}_{2,\varrho}(Q_{T})$. Set
\[
C^{n}_{s,t}=\int_{s}^{t}(f_n^{0}
+\mbox{div}\bar{f}_n)(\theta,X_\theta)\,d\theta,\quad
R^{n}_{s,t}=\int_{s}^{t}(g_n^{0}
+\mbox{div}\bar{g}_n)(\theta,X_\theta)\,d\theta,
\]
for $0\le s\le t\le T$ and
\begin{align}\label{eqs.n1}
C_{r,t}&=\int_{r}^{t}f^{0}(\theta,X_{\theta})\,d\theta
+\int_{r}^{t}a^{-1}\bar{f}(\theta,X_\theta)\,d^{*}X_{\theta},\\
\label{eqs.n2}
R_{r,t}&=\int_{r}^{t}g^{0}(\theta,X_{\theta})\,d\theta
+\int_{r}^{t}a^{-1}\bar{g}(\theta,X_\theta)\,d^{*}X_{\theta}
\end{align}
for $0\le s<r\le t\le T$. It is clear that  for every $(s,x)\in
Q_{\hat{T}}$\,,
\[
A^{u}_{r,t}=C_{r,t}-R_{r,t},\quad 0\le s<r\le t\le T,\quad
P_{s,x}\mbox{-}a.s..
\]
By (\ref{inl.1}), (\ref{inl.2}), for every $(s,x)\in Q_{\hat{T}}$
and $r\in(s,T]$,
\begin{equation}
\label{eq.eq} E_{s,x}\sup_{r\le t\le T} (|C^{n}_{r,t}-C_{r,t}|
+|R^{n}_{r,t}-R_{r,t}|) \rightarrow0,
\end{equation}
which implies (\ref{eq6.1}). Now, let $v\in\mathcal{W}_{\varrho}$
be such that  $\mathcal{L}v=\mu_{1}$ and $v(T)=0$. By
(\ref{eq6.1}), there exists CAF $\tilde{C}$ such that
$X^{v}=M^{v}+\tilde{C}$ in the sense of (\ref{eq6.1}). Since
$C,\tilde{C}$ satisfy (\ref{eqs.n1}), $C=\tilde{C}$. Hence, by
Aronson's upper estimate and a priori estimates for PDEs,
\begin{align*}
E_{s,x}|C_{r,T}|^2&\le C(E_{s,x}|M^{v}_{r,T}|^{2}
+E_{s,x}|v(r,X_{r})|^2\\
&\le C\frac{\varrho^{2}(x)}{(r-s)^{d/2}}(\|\nabla
v\|^{2}_{2,\varrho,T} +\sup_{0\le t\le T}\|v(t)\|^{2}_{2,\varrho})
\le C\frac{\varrho^{2}(x)}{(r-s)^{d/2}}\|\mu_{1}\|_{*}\,,
\end{align*}
which proves (\ref{eq6.1.5}). To show (\ref{eq6.2}), (\ref{eq6.3})
let us fix $(s,x)\in Q_{\hat{T}}$, $r\in(s,T]$  and choose $\xi\in
C^\infty_{0}(Q_T)$  so that $\xi{\mathbf{1}}_{Q_{s+\delta}}=0$ for
some $\delta\in (0,T-s)$. Then, by Proposition \ref{prop.p},
$\eta=\xi p(s,x,\cdot,\cdot)\in
{\mathbb{L}}_2(0,T;H^{1}_{\varrho})$ and
\begin{align*}
&E_{s,x}\int_{s+\delta}^T\xi(t,X_t)\,dC_{s,t}
=\lim_{n\rightarrow\infty}
E_{s,x}\int_{s+\delta}^T\xi(t,X_t)\,dC^{n}_{s,t}\\
&\quad= \lim_{n\rightarrow\infty}\Phi_1^n(\eta) =\Phi_1(\eta)
=\int_{s+\delta}^T\int_{{\mathbb{R}}^d} \xi(t,y)p(s,x,t,y)\,d\mu_1
\end{align*}
by (\ref{eq6.1.5}), (\ref{eq.eq}) and the fact that
$\mu_1^{n}\rightarrow \mu_1$ in
${\mathbb{L}}_{2}(0,T;H_{\varrho}^{-1})$. From this we easily get
(\ref{eq6.2})  and (\ref{eq6.3}). Passing to the limit with
$r\rightarrow s^{+}$ in (\ref{eq6.1}) and using  (\ref{eqs.n1}),
(\ref{eqs.n2}) and Corollary \ref{coll.1} we get (\ref{eq6.1})
with $r=s$ for q.e. $(s,x)\in Q_{\hat{T}}$. Similarly, passing to
the limit with $r\rightarrow s^{+}$ in (\ref{eq6.2}) and
(\ref{eq6.3}) and using (\ref{RM1}) we get (\ref{eq6.2}) and
(\ref{eq6.3}) with $r=s$ for q.e. $(s,x)\in Q_{\hat{T}}$.
\end{proof}
\medskip

From now on we write $C\sim \mu$ if CAF $C$ is associated with
measure $\mu$ in the sense of (\ref{eq6.2}). From the above
theorems we get in particular the well known Revuz correspondence
for smooth measures. However in the case of diffusion process
$(X,P_{s,x})$ this correspondence might be expressed via density
of the process which we present in the following corollary.
\begin{remark}\label{rems.1}
{\rm Repeating  proofs of Lemmas 2.2.8 and 2.2.9 in
\cite{Fukushima} with  one can show that if $\mu\in S^{c}$ then
there exists a sequence $\{F_{n}\}$ (called nest) of closed
subsets of $\check{Q}_{T}$ such that
$\mu(\check{Q}_{T}\setminus\bigcup_{n=1}^{+\infty}F_{n})=0$,
$\lim_{n\rightarrow\infty}\overline{\mbox{\rm cap}}(K-F_{n})=0$
for every compact $K\subset\check{Q}_{T}$ and
$\mathbf{1}_{F_{n}}d\mu\in S^{c}_{0}$ for every $n\in\mathbb{N}$.
}
\end{remark}

\begin{definition}
{\rm We say that $dK:\Omega\times \mathcal{B}([0,T])\rightarrow
\mathbb{R}$ is a random measure if
\begin{enumerate}
\item[(a)]$dK(\omega)$ is a measure
for every $\omega\in\Omega$,
\item[(b)] $\omega\mapsto dK(\omega)$
is $(\mathcal{G},\mathcal{B}(\mathcal{M}[0,T]))$-measurable,
\item[(c)]$\int_{s}^{t}\,dK_{\theta}$ is
$\mathcal{G}^{s}_{t}$-measurable for every $0\le s\le t\le T$.
\end{enumerate}
}
\end{definition}

\begin{remark}
{\rm By results proved in \cite{Y.Osh.DirichletForms} one can
associate with  the operator $\mathcal{L}$ a Hunt process
$\{(Z_{t},\tilde{P}_{z}),t\ge 0$, $z\in\mathbb{R}^{d+1}\}$.
Actually, it follows from \cite{Y.Osh.DirichletForms} that
$\tilde{P}_{z}$ coincides with $P_{s,x}$ for $z=(s,x)\in
Q_{\hat{T}}$ and that $Z_{t}=(\tau(t), X_{\tau(t)})$, where $\tau$
is the uniform motion to the right, i.e. $\tau(t)=\tau(0)+t$ and
$\tau(0)=s$ under $P_{s,x}$. }
\end{remark}

\begin{lemma}\label{lmd.3}
Let $\{dK^{n}\}$ be a sequence of random measures. Assume that for
$(s,x)\in F\subset Q_{\hat{T}}$ there exist random elements
$dK^{s,x}:(\Omega,\mathcal{G})\rightarrow (\mathcal{M}^{+}[0,T],
\mathcal{B}(\mathcal{M}^{+}[0,T])$such that
\[
dK^{n}(\cdot,X_{\cdot})\rightarrow dK^{s,x} \mbox{ in
}\mathcal{M}^{+}([0,T]\mbox{ in probability } P_{s,x}.
\]
Then there exists a random measure $dK$ such that
\[
dK^{s,x}=dK,\quad P_{s,x}\mbox{-}a.s.
\]
for every $(s,x)\in F$.
\end{lemma}
\begin{proof}
Let $n_{0}(s,x)=0$ and let
\[
n_{k}(s,x)=\inf\{m>n_{k-1}(s,x), \sup_{p,q\ge
m}P_{s,x}(d_{M}(dK^{p},dK^{q})>2^{-k})<2^{-k}\}
\]
dor $k\ge1$. By induction, $n_{k}\in\mathcal{B}(Q_{\hat{T}})$ for
every $k\ge0$, and hence  $dL^{s,x,k}=dK^{n_{k}(s,x)}$ is
$\mathcal{B}(Q_{\hat{T}})\otimes
\mathcal{G}/\mathcal{B}(\mathcal{M}^{+}[0,T])$ measurable. Put
\begin{equation}
\label{eq1.3} dL^{s,x}(\omega)= \left\{
\begin{array}{ll} \lim_{k\rightarrow\infty}
dL^{s,x,k}(\omega)\mbox{ in } \mathcal{M}([0,T]),&
\mbox{ if the limit exists},\\
0,&\mbox{otherwise}.
\end{array}
\right.
\end{equation}
By the Borel-Cantelli lemma, for every $(s,x)\in F$ the limit in
(\ref{eq1.3}) exists $P_{s,x}$-a.s. and $dL^{s,x}=dK^{s,x}$,
$P_{s,x}$-a.s.. Putting $dK(\omega)=dL^{Z_{0}(\omega)}$ we get
random measure having the desired properties.
\end{proof}

Let $\mu\in S^{c}$. In what follows by $d\mu(\cdot,X_{\cdot})$ we
denote random measure such that for q.e. $(s,x)\in Q_{\hat{T}}$,
\begin{eqnarray}\label{R}
E_{s,x}\int_s^T\xi(t,X_t)\,d\mu(t,X_{t})
=\int_{Q_{sT}}\xi(t,y)p(s,x,t,y)\,d\mu(t,y)
\end{eqnarray}
for every $\xi \in \mathcal{B}^{+}(Q_{T})$.

\begin{corollary}\label{m}
For every $\mu\in S^{c}$ there exists a unique random measure
$d\mu(\cdot,X_{\cdot})$. Moreover, for every $\mu\in S^{c}_{0}$
and $s\in[0,T)$,
\begin{equation}
\label{eq5.7}
E_{s,\varrho}\left(\int_{s}^{T}\,d\mu(t,X_{t})\right)^2\le
C\|\mu\|^{2}_{*}\,.
\end{equation}
\end{corollary}
\begin{proof}
Uniqueness follows from Proposition \ref{LT} and Remark
\ref{rems.1}. Let  $\mu\in S^{c}_{0}$ and let $u$ be a unique
solution of PDE$(0,\mu)$. By Proposition  \ref{tw3.1} there exist
a unique positive CAF $A^{\mu}$ such that $A^{\mu}\sim\mu$ and a
version of $u$ (still denoted by $u$) such that for every
$(s,x)\in Q_{\hat{T}}$ and $r\in (s,T]$,
\begin{equation}
\label{eqA.2} u(r,X_{r})=A^{\mu}_{r,T}-M^{u}_{r,T}, \quad 0\le
s<r\le T,\quad P_{s,x}\mbox{-}a.s..
\end{equation}
By Theorem \ref{tw3.1}, the  random measures
\begin{equation*}
dK^{n}(\omega)= \left\{
\begin{array}{ll} \Phi_{n}(t,X_t(\omega))\,dt,&\mbox{if }
\int_{0}^{T}\Phi_{n}(t,X_t(\omega))\,dt<\infty,\\
0,& \mbox{otherwise},
\end{array}
\right.
\end{equation*}
where $\Phi_{n}$ are defined as in the proof Theorem \ref{tw3.1},
satisfy the assumptions of Lemma \ref{lmd.3}. Hence there exists a
unique random measure $d\mu(\cdot,X_{\cdot})$ such that
$\int_{s}^{t}d\mu(\theta,X_{\theta})=A^{\mu}_{s,t},\,\, s\le t\le
T,\,\, P_{s,x}$-a.s. for q.e. $(s,x)\in Q_{\hat{T}}$. Therefore,
by (\ref{eqA.2}),
\begin{equation}
\label{eqA.1}
u(r,X_{r})=\int_{r}^{T}d\mu(\theta,X_{\theta})-M^{u}_{r,T}, \quad
0\le s<r\le T,\quad P_{s,x}\mbox{-}a.s..
\end{equation}
Integrating (\ref{eqA.1}) with respect to $\varrho^{2}\,dm$ and
using  Proposition \ref{prop3.1} yields
\[
E_{s,\varrho}(\int_{r}^{T}d\mu(\theta,X_{\theta}))^{2} \le
C(\|u(r)\|^{2}_{2,\varrho,T} +\|\nabla u\|^{2}_{2,\varrho,T})\le
C\|\mu\|_{*}^{2}
\]
for every  $r\in (s,T]$, the last inequality being a consequence
of Theorem \ref{tha.1}. The result now follows from Fatou's lemma.
Now, let $\mu\in S^{c}$. Then, by Remark \ref{rems.1}, there
exists a nest $\{F_{n}\}$ such that
$\mu_{n}=\mathbf{1}_{F_{n}}d\mu\in S^{c}_{0}$. By what has already
been proved, for each $n\in{\mathbb{N}}$ there exists the random
measure $d\mu_{n}(\cdot,X_{\cdot})$. Let us observe that if $n\le
m$ then $\mathbf{1}_{F_{n}}d\mu_{m}=d\mu_{n}$, which  implies that
$\mathbf{1}_{F_{n}}d\mu_{m}(\cdot,X_{\cdot})=d\mu_{n}(\cdot,X_{\cdot})$.
Therefore, $d\mu_{m}(\cdot,X_{\cdot})\ge
d\mu_{n}(\cdot,X_{\cdot})$. By Lemma \ref{lmd.3} it follows that
there exists random measure $dK$ such that
$dK=\lim_{n\rightarrow\infty}d\mu_{n}(\cdot,X_{\cdot})$ in
$\mathcal{M}^{+}[0,T]$ in probability $P_{s,x}$ for q.e. $(s,x)\in
Q_{\hat{T}}$. It is clear that $dK$ satisfies (\ref{R}). Therefore
$dK=d\mu(\cdot,X_{\cdot})$.
\end{proof}

\begin{remark}\label{rems.2}
{\rm Let $u$ satisfies the assumptions of Proposition \ref{tw3.2}.
Then by (\ref{eq6.1}) and a priori estimates for BSDEs (see
\cite{EKPPQ}), for every $(s,x)\in Q_{\hat{T}}$ and $r\in (s,T]$,
\begin{align*}
&E_{s,x}\sup_{r\le t\le T}|u(t,X_{t})|^2
+E_{s,x}\int_{r}^{T}|\nabla u|^{2}(\theta,X_{\theta})\,d\theta\\
&\qquad \le C\left(E_{s,x}|u(T,X_{T})|^2
+E_{s,x}(\int_{r}^{T}\,d\mu_{1}(\theta,X_{\theta}))^{2}
+E_{s,x}(\int_{r}^{T}\,d\mu_{2}(\theta,X_{\theta}))^{2}\right).
\end{align*}
Hence, if $E_{s,x}(\int_{r}^{T}d\mu_{1}(\theta,X_{\theta}))^{2}
+E_{s,x}(\int_{r}^{T}d\mu_{2}(\theta,X_{\theta}))^{2}<\infty$,
then (\ref{eq6.1}), (\ref{eq6.2}), (\ref{eq6.3}) are satisfied
with $r=s$. Consequently, by Corollary \ref{m}, for each fixed
$s\in [0,T)$, (\ref{eq6.1}), (\ref{eq6.2}), (\ref{eq6.3}) are
satisfied for a.e. $x\in\mathbb{R}^{d}$. If $s\in(0,T)$ this
follows also from the fact that cap$_{L}(\{s\}\times B)>0$ for
every $B\in{\mathcal{B}}(Q_T)$ such that  $m(B)>0$. }
\end{remark}

\begin{definition}
{\rm We say that $X^u$ is a locally finite semimartingale if it is
a semimartingale under $P_{s,x}$ for q.e. $(s,x)\in Q_{T}$ and its
finite variation part is a locally finite CAF.}
\end{definition}

Let us remark that the  class of locally finite semimartingales
appears naturally when considering Revuz duality for additive
functionals (see \cite{Fukushima1}).

The next theorem shows that the condition
$\mathcal{L}u\in\mathcal{M}$ is necessary and sufficient for $X^u$
to be locally finite semimartingale.

\begin{theorem}
Let $u\in \mathcal{W}_{\varrho}$.
\begin{enumerate}
\item[\rm(i)]$\mathcal{L}u\in\mathcal{M}$ iff $X^{u}$ is a locally finite
semimartingale.
\item[\rm(ii)]Assume that $\mathcal{L}u\in\mathcal{M}$. Let $\mu=\mathcal{L}u$
and let $A^{u}$ denote the finite variation part  of $X^u$. Then
$dA^{u}=d\mu(\cdot,X_{\cdot})$.
\end{enumerate}
\end{theorem}
\begin{proof}
Suppose that $\mathcal{L}u\in\mathcal{M}$ and let $\mu=\mathcal{L}u$. From Theorem \ref{decom1} it follows that
$\mu\ll \overline{{\rm cap}}$. Let $\mu=\mu^+-\mu^-$ be the
canonical decomposition. Of
course, $\mu^{+}\mu^{-} \ll \overline{\mbox{cap}}$. Hence, by
Theorem \ref{decom}, there exist $\gamma_{1},\gamma_{2}\in S^{c}_{0}$ and
$\alpha_{1},\alpha_{2}\in {\mathbb{L}}^{+}_{1,loc}(Q_{T})$ such
that $\mu^{+}=\alpha_{1}\,d\gamma_{1}$,
$\mu^{-}=\alpha_{2}\,d\gamma_{2}$. Put
\[
D^{u}_{s,t}=\int_{s}^{t}\alpha_{1}
(\theta,X_{\theta})\,d\gamma_{1}(\theta,X_{\theta})
-\int_{s}^{t}\alpha_{2}(\theta,X_{\theta})\,d\gamma_{2}
(\theta,X_{\theta}),\quad
0\le s\le t\le T.
\]
By Aronson estimates,  for every $\eta\in C_{0}(Q_{T})$,
\[
\int_{Q_{T}}(E_{s,x}\int_{s}^{T}
\eta\alpha_{i}(t,X_t)\,d\gamma_{i}(t,X_t)) \le C
\int_{Q_{T}}\eta\alpha_{i}\,d\gamma_{i},\quad i=1,2.
\]
From the above and Proposition \ref{propl.1} it follows that for
q.e. $(s,x)\in Q_{\hat{T}}$ the functional $D^{u}$ is well defined
and $E_{s,x}\int_{s}^{T}
\eta\alpha_{i}(t,X_t)\,d\gamma_{i}(t,X_t)<\infty$ for every
$\eta\in C_{0}(Q_{T})$. By Theorem \ref{tw3.1},
\[
X^{u}_{s,t}=M^{u}_{s,t}+A^{u}_{s,t},\quad t\in[s,T],\quad
P_{s,x}\mbox{-}a.s.
\]
for q.e. $(s,x)\in Q_{\hat{T}}$ with $A^{u},M^{u}$ as in Theorem
\ref{tw3.1}.  We shall show that $A^{u}=D^{u}$. In view of
Proposition \ref{LT}, to prove this it suffices to show that for
every $\eta\in {\mathcal{W}}_{\varrho}\cap C_{c}(Q_{T})$,
\begin{eqnarray}\label{LL}
E_{s,x}\int_{s}^{T}\eta(t,X_t)\,dA^{u}_{s,t}
=E_{s,x}\int_{s}^{T}\eta(t,X_t)\,dD^{u}_{s,t}
\end{eqnarray}
for $(s,x)\in D_{0}(\eta\cdot A^{u})\cap D_{0}(\eta\cdot D^{u})$.
Given  $\delta\in(0,T-s)$ write
\begin{eqnarray}\label{eq.delta}
L_{\delta}=E_{s,x}\int_{s+\delta}^{T}\eta(t,X_t)\,dA^{u}_{s,t},\quad
R_{\delta}=E_{s,x}\int_{s+\delta}^{T}\eta(t,X_t)\,dD^{u}_{s,t}.
\end{eqnarray}
Then by Theorem \ref{tw3.1},
\begin{align*}
L_{\delta}&=\int_{Q_{s+\delta,T}}\eta f^{0}(t,y)p(s,x,t,y)\,dt\,dy
-\int_{Q_{s+\delta,T}}\nabla\eta \bar{f}(t,y)p(s,x,t,y)\,dt\,dy\\
&\quad+
\int_{Q_{s+\delta,T}}\eta\bar{f}(t,y)\nabla_{y}p(s,x,t,y)\,dt\,dy
\end{align*}
and
\[
R_{\delta}=\int_{Q_{s+\delta,T}}\alpha_{1}
\eta(t,y)p(s,x,t,y)\,d\gamma_{1}(t,y)
-\int_{Q_{s+\delta,T}}\alpha_{2}\eta(t,y)p(s,x,t,y)\,d\gamma_{2}(t,y).
\]
Since $\eta p(s,x,\cdot,\cdot)\in
C_{0}(Q_{s+\delta,T})\cap{\mathbb{L}}_{2}(s+\delta,T;H^{1}_{\varrho})$,
it follows from the assumption that $L_{\delta}=R_{\delta}$ for
every $\delta>0$. Letting $\delta\rightarrow 0^{+}$ for $(s,x)\in
D_{0}(\eta\cdot A^{u})\cap D_{0}(\eta\cdot D^{u})$ we get
$L_{0}=R_{0}$.

Now,  assume that $X^{u}$ is a locally finite semimartingale.
Without lost of generality we may and will assume that $b=0$. The
general case can be handled easily by using Girsanov's theorem,
because under the change of measure removing the drift term in the
decomposition of $X^u$  new terms of finite variation appear (see
e.g. \cite[Section 4]{Roz.Div} for details). Then $A^{u}$ from the
decomposition of $X^{u}$ of Theorem \ref{tw3.1} is of finite
variation. Given $\eta\in C_{c}(Q_{T})$ put
\[
\mu(\eta)=\int_{{\mathbb{R}}^{d}}(E_{0,x}\int_{0}^{T}
\eta(\theta,X_{\theta})\,dA_{0,\theta}).
\]
By the assumption, the above integral is well defined and the
functional $\mu$ is continuous with respect to the uniform
convergence on compacts, which implies that $\mu$ is a measure. We
shall show that $\mathcal{L}u=\mu$. Let $\eta\in
R_{\alpha}(C_{c}^{\infty}(Q_{T})) \subset
D(\mathcal{L})\subset\mathcal{W}_{\varrho}$, where $R_{\alpha}$ is
the resolvent of  $\mathcal{L}$. Then, by Theorem \ref{tw3.1},
\[
\eta(t,X_{t})=\eta(T,X_{T})-\int_{t}^{T}{\mathcal{L}}
\eta(\theta,X_{\theta})\,d\theta-\int_{t}^{T}
\sigma\nabla\eta(\theta,X_{\theta})\,dB_{s,\theta},\quad t\in
[s,T]
\]
$P_{s,x}$-a.s. for q.e. $(s,x)\in Q_{\hat{T}}$.  Integrating by
parts we get
\begin{align}
\label{eq6.13} E_{0,x}u(0,X_{0})\eta(t,X_{t})
&=E_{0,x}u(T,X_{T})\eta(T,X_{T})-E_{0,x}\int_{0}^{T}
u(t,X_t){\mathcal{L}}\eta(t,X_t)\,dt \nonumber\\
&\quad -E_{0,x}\int_{0}^{T}\eta(t,X_t)\,dA_{0,t}
-E_{0,x}\int_{0}^{T}\langle a\nabla \eta,\nabla
u\rangle(t,X_t)\,dt.
\end{align}
Notice that $R_{\alpha}\xi\varrho^{-2}\in\mathbb{L}_{2}(Q_{T})$ if
$\xi\in C_{c}^{\infty}(Q_{T})$. This follows from Proposition
\ref{prop3.1} and the fact that
\[
R_{\alpha}\xi(s,x)=E_{s,x}\int_{0}^{T}
\mathbf{1}_{[0,T]}(s+t)e^{-\alpha t}\xi(s+t,X_{s+t})\,dt
\]
(see, e.g., \cite{Oshima1}). Integrating (\ref{eq6.13}) with
respect to $x$ and using symmetry of the operator $L_t$ we get
\[
\langle u(0),\eta(0)\rangle_{2} =\langle u(T),\eta(T)\rangle_{2}
-\langle u,\frac{\partial \eta}{\partial t}\rangle_{2,T} +\langle
u, L_{t}\eta\rangle_{2,T}-\int_{Q_{T}}\eta\,d\mu,
\]
which proves that $\langle\mathcal{L}u,\eta\rangle_{2,T}
=\int_{Q_{T}}\eta\,d\mu$ for all $\eta\in
R_{\alpha}(C^{\infty}_c(Q_T))$. That $\mathcal{L}u=\mu$ now
follows from strong continuity of the resolvent.
\end{proof}

Using Theorem \ref{tw3.2} one can prove useful estimate for the
first moment of the supremum of $X^{u}$ in terms of the norm of
$u$ in ${\mathcal{W}}_{\varrho}$.
\begin{corollary}\label{f}
If $u\in {\mathcal{W}}_{\varrho}$ then there is a
quasi-continuous version of $u$ (still denoted by $u$) such that for every $s\in (0,T)$
\[
E_{s,\sqrt{\varrho}}\sup_{s\le t\le T}|u(t,X_{t})| \le
C\|u\|_{{\mathcal{W}_{\varrho}}}.
\]
\end{corollary}
\begin{proof}
Since $\partial u/\partial
t\in{\mathbb{L}}_2(0,T;H^{-1}_{\varrho})$, $u$
admits the representation
\begin{equation}
\label{eq3.9}
\frac{\partial u}{\partial t}+L_{t}u
=f^{0}+\mbox{div}(\bar{f})+\frac12\mbox{div}(a \nabla u)+b\nabla u
\end{equation}
for some $f^{0},\bar{f}\in{\mathbb{L}}_{2,\varrho}(Q_{T})$. From
Theorem \ref{tw3.1} it follows that there exists a
quasi-continuous  version of $h$ (still denoted by $h$) such that
\begin{align*}
X^u_{s,t}&=M^u_{s,t}+\int_{s}^{t}f^{0}(\theta,X_{\theta})\,d\theta
+\int_{s}^{t}a^{-1}\bar{f}(\theta,X_{\theta})\,d^{*}X_{\theta}\\
&\quad+\frac12\int_{s}^{t}\nabla
u(\theta,X_{\theta})\,d^{*}X_{\theta} +\int_{s}^{t}\nabla
u(\theta,X_{\theta})d\beta_{s,\theta} \nonumber
\end{align*}
for q.e. $(s,x)\in Q_{\hat{T}}$. By Doob's
${\mathbb{L}}_2$-inequality,
\begin{align*}
E_{s,x}\sup_{s\le t\le T}|u(t,X_{t})| &\le
C\left(E_{s,x}(|u(T,X_T)|^2 +\int_s^T(|f^0|^2
+|\nabla u|^2)(\theta,X_\theta)\,d\theta)\right)^{1/2}\\
&\quad+\sup_{s\le t\le T} \left|\int_{t}^{T} (a^{-1}\bar{f}+\nabla
u) (\theta,X_{\theta})\,d^{*}X_{\theta}\right|.
\end{align*}
Multiplying the above inequality by $\varrho$ and using
Proposition \ref{prop3.2} we obtain
\[
E_{s,\sqrt{\varrho}}\sup_{s\le t\le T}|u(t,X_{t})|
\le C(\|u(T)\|_{2,\varrho}+\|f^{0}\|_{2,\varrho,T}
+\|\bar{f}\|_{2,\varrho,T}+\|\nabla u\|_{2,\varrho,T}).
\]
Taking infimum over all
$f^{0},\bar{f}\in{\mathbb{L}}_{2,\varrho}(Q_{T})$ such that
(\ref{eq3.9}) is satisfied yields
\[
E_{s,\sqrt{\varrho}}\sup_{s\le t\le T}|u(t,X_{t})|\le
C(\|u(T)\|_{2,\varrho} +\|\frac{\partial u}{\partial
t}\|_{*}
+\|\nabla u\|_{2,\varrho,T}).
\]
This proves the desired estimate because the imbedding of
${\mathcal{W}}_{\varrho}$ into the vector space
$C([0,T],{\mathbb{L}}_{2,\varrho}({\mathbb{R}}^{d}))$ is
continuous (see, e.g., \cite{Lions}).
\end{proof}

To estimate the second moment of the supremum of $X^{u}$ we assume
that $\mathcal{L}u\in S^{c}_{0}-S^{c}_{0}$. It is worth noting
that solutions  of parabolic equations with the right-hand side in
${\mathbb{L}}_{2,\varrho}(Q_{T})$ and solutions of unilateral or
bilateral problems satisfy that assumption.
\begin{corollary}
Let $u\in {\mathcal{W}}_{\varrho}$ and $\mathcal{L}u\in
S^{c}_{0}-S^{c}_{0}$. Then there is a quasi-continuous version of
$u$ (still denoted by $u$) such that for every $s\in (0,T)$,
\[
E_{s,\varrho}\sup_{s\le t\le T}|u(t,X_{t})|^{2} \le
C(\|\mu^{+}\|^{2}_{*}
+\|\mu^{-}\|^{2}_{*}),
\]
where
$\mu^+,\mu^-\in S^{c}_{0}$ and
$\mathcal{L}u=\mu^+-\mu^-$.
\end{corollary}
\begin{proof}
By Theorem \ref{tw3.2}, $X^u$ admits the decomposition
(\ref{eq6.1})  for q.e. $(s,x)\in Q_{\hat T}$. Therefore one can
prove the desired estimate by the same method as in the proof of
Corollary \ref{f}.
\end{proof}

\nsubsection{Appendix}
\label{sec7}

For convenience of the reader we collect here some estimates for
diffusions ${\mathbb{X}}$ associated with $L_t$ and related
estimates on the  fundamental solution $p$ of $L_t$ and  weak
solutions of the Cauchy problem
\begin{equation}
\label{eq4.6} \frac{\partial u}{\partial t}+L_tu
=-\Phi,\quad u(T)=\varphi
\end{equation}
(PDE$(\varphi,\Phi)$ for short), where
$\Phi\in\mathbb{L}_{2}(0,T;H^{-1}_{\varrho})$.
Recall that
$u\in\mathcal{W}_{\varrho}$ is a strong
solution of PDE$(\varphi,\Phi)$ if for any
$\eta\in\mathbb{L}_{2}(0,T;H^{1}_{\varrho})$,
\begin{align*}
&\int_{t}^{T}\langle \frac{\partial u}{\partial s}(s),
\eta(s)\rangle_{\varrho}\,ds
-\frac12\int_{t}^{T}\langle a(s)\nabla
u(s),\nabla (\varrho^{2}\eta(s))\rangle_{2}\,ds
+\int_{t}^{T}\langle b(s)\nabla u(s),\eta(s)\rangle_{2,\varrho,T}\\
&\qquad=\int_{t}^{T}\langle f^0(s),\eta(s)\rangle_{2,\varrho}\,ds
-\int_{t}^{T}\langle
\bar{f}(s),\nabla(\varrho^{2}\eta(s))\rangle_{2}\,ds .
\end{align*}
for all $t\in[0,T]$, where
$f^{0},\bar{f}\in\mathbb{L}_{2,\varrho}(Q_{T})$ are such  that
$\Phi=f^{0}+\dyw\bar{f}$.

For the proof of the following theorem  see e.g. \cite{Lad,Lions}.
\begin{theorem}\label{tha.1}
For every $\Phi\in\mathbb{L}_{2}(0,T,H^{-1}_{\varrho})$ there
exists a unique  strong solution $u\in\mathcal{W}_{\varrho}$ of
PDE$(\varphi,\Phi)$ and
\[
\sup_{0\le t\le T}\|u(t)\|^{2}_{2,\varrho}+\|\nabla
u\|^{2}_{2,\varrho,T} \le
C(\|\varphi\|^{2}_{2,\varrho}+\|\Phi\|^{2}_{*}).
\]
\end{theorem}

\begin{proposition}
\label{prop3.1} Let $\varrho\in\mathcal{R}$. There exist $0<C_1\le
C_2$  depending only on $\lambda,\Lambda,\Lambda_{1},d,T$ and
$\varrho$ such that
\begin{align*}
C_1\int_t^T\!\!\int_{{\mathbb{R}}^d}
|\psi(\theta,x)|\varrho(x)\,d\theta\,dx
&\le\int_t^T\!\!\int_{{\mathbb{R}}^d}
E_{s,x}|\psi(\theta,X_{\theta})|\varrho(x)\,d\theta\,dx \\
& \le C_2\int_t^T\!\!\int_{{\mathbb{R}}^d}
|\psi(\theta,x)|\varrho(x)\,d\theta\,dx
\end{align*}
for any $\psi\in{\mathbb{L}}_{1,\varrho}(Q_T)$ and $t\in[s,T]$.
\end{proposition}
\begin{proof}
Follows from Proposition 5.1 in Appendix in \cite{V.Bally} and
Aronson's estimates (see \cite[Theorem 7]{Aronson}).
\end{proof}
\medskip

We now provide  useful estimates for  moments of ${\mathbb{X}}$.

\begin{lemma}
\label{moments} For every $p\ge1$ there is $C$ depending only on
$\lambda,\Lambda,d,T$ and $p$ such that
\[
E_{s,x}\sup_{s\le t\le T}|X_t|^p\le C(p)(1+|x|)^{p}.
\]
\end{lemma}
\begin{proof}
By \cite{Stroock.DIFF} there exist $C_1,C_2>0$ such that for every
$(s,x)\in Q_{\hat{T}}$ and $r\ge 0$,
\[
P_{s,x}(\sup_{s\le t\le T}|X_{t}-x|> r)\le
C_{1}\exp(\frac{-C_{2}r^{2}}{T-s}).
\]
From this we conclude that for every $p\ge 0$,
\[
E_{s,x}\sup_{s\le t\le T}|X_{t}-x|^{p}\le C(p)
\]
from which the result follows.
\end{proof}
\medskip

The following estimates for $p$ and weak solutions of
(\ref{eq4.6}) are known, but originally stated in terms of
${\mathbb{L}}_{p,q,\varrho}$\,-norms with $\varrho\equiv1$. At the
expense of minor technical changes their proofs my be adapted to
the case of spaces with weight $\varrho$ such that $\varrho^{-1}$
is a polynomial. For the first proposition see Theorems 5, 7 and
10, and for the second one Theorems 5 and 10 in \cite{Aronson}.

\begin{proposition}\label{prop.p}
Assume that $p, q\in(1,+\infty]$, $\frac{d}{2p}+\frac{1}{q}<1$.
Then for any $(s,x)\in Q_{\hat T}$,
\begin{enumerate}
\item [\rm(i)] $\|p(s,x,\cdot,\cdot)\|_{p^{'},q^{'},\varrho^{-1}_{x}}
+\|\nabla p(s,x,\cdot,\cdot)\|_{(2p)^{'},(2q)^{'},\varrho^{-1}_{x}}<C$,
\item [\rm(ii)]  $p(s,x,\cdot,\cdot),\nabla p(s,x,\cdot,\cdot)
\in {\mathbb{L}}_{2,\varrho^{-1}}(Q_{s+\delta,T})$
for every $\delta\in (0,T-s]$.
\end{enumerate}
\end{proposition}

\begin{proposition}\label{prop.pp}
Let $p,q$ satisfy the assumption of Proposition \ref{prop.p}. Then
there exists a continuous version $u$ of a weak solution of
(\ref{eq4.6}), and
\[
|u(t,x)|\le
C\varrho^{-1}(x)(\|\varphi\|_{\infty,\varrho}+\|f\|_{p,q,\varrho}
+\|\bar{f}\|_{2p,2q,\varrho}),\quad (t,x)\in Q_{\hat T}.
\]
\end{proposition}

\begin{proposition}\label{prop3.2}
Assume that $p, q\in(1,\infty]$, $\frac{d}{2p}+\frac1q<\frac12$\,.
\begin{enumerate}
\item[\rm(i)] For every $(s,x)\in Q_{\hat T}$
and $f\in{\mathbb{L}}_{p,q,\varrho}(Q_T)$,
\[
E_{s,x}\int_{s}^{T}|f(t,X_t)|^2\,dt
\le C\varrho^{-2}(x)\|f\|_{p,q,\varrho}\,.
\]
\item [\rm(ii)] For every $(s,x)\in
Q_{\hat T}$ and $\bar f\in{\mathbb{L}}_{2p,2q,\varrho}(Q_T)$ the
integral (\ref{dec1}) is well defined and
\[
E_{s,x}\sup_{s\le t\le T}
|\int_{s}^{t}\bar{f}(\theta,X_{\theta})\,d^{*}X_{\theta}|^{2} \le
C \varrho^{-2}(x)\|\bar{f}\|_{2p,2q,\varrho}\,.
\]
\item[\rm(iii)] For every $s\in [0,T)$ and
$f\in{\mathbb{L}}_{2,\varrho}(Q_T)$,
\[
E_{s,\varrho}\int_{s}^{T}|f(t,X_t)|^2\,dt\le
C\|f\|_{2,\varrho,T}\,.
\]
\item[\rm(iv)] If $\bar{f}\in{\mathbb{L}}_{2,\varrho}(Q_T)$ then
te integral (\ref{dec1}) is well defined for a.e. $(s,x)\in
Q_{\hat{T}}$\,,
\begin{equation}
\label{eq4.3} \int_{Q_{T}}\left(E_{s,x}\int_{s}^{T}
|\bar{f}(t,X_t)|\,
d|\alpha^{s,x}_{s,\cdot}|_{\theta}\right)\varrho(x)\,dx\,ds \le
C\|\bar{f}\|_{2,\varrho,T}
\end{equation}
and
\[
\int_{Q_{T}}(E_{s,x}\sup_{s\le t\le T}|\int_{s}^{t}
\bar{f}(\theta,X_{\theta})\,d^{*}X_{\theta}|)\varrho(x)\,dx\,ds
\le C\|\bar{f}\|_{2,\varrho,T}
\]
\end{enumerate}
\end{proposition}
\begin{proof}
(i) Since $\varrho^{-1}(x+y)\le C\varrho^{-1}(x)\varrho^{-1}(y)$,
applying H\"older's inequality gives
\begin{align*}
\int_{Q_{sT}}|f(t,y)|^2p(s,x,t,y)\,dt\,dy &\le
C\varrho^{-2}(x)\int_{Q_{sT}}|f(t,y)|^2
\varrho^2(y) p(s,x,t,y)\varrho_{x}^{-2}(y)\,dt\,dy\\
&\le C\varrho^{-2}(x)\|f\|^{2}_{p,q,\varrho}
\|p(s,x,\cdot,\cdot)\varrho_{x}^{-2}\|_{(p/2)^{'},(q/2)^{'}},
\end{align*}
and the result follows from Proposition \ref{prop.p}.
\\
(ii) By (i) integrals with respect to backward and forward
martingale  are well defined. As for the finite variation part
observe that
\begin{align*}
E_{s,x}\int_{s}^{T}
|\bar{f}|(\theta,X_{\theta})\,d|\alpha^{s,x}_{s,\cdot}|_{\theta}
&=\int_{Q_{s,T}}|\bar{f}(\theta,y)| |\nabla
p|(s,x,\theta,y)\,d\theta\,dy\\
& \le
C\varrho^{-1}(x)\|\bar{f}\|_{2p,2q,\varrho} \|\nabla
p(s,x,\cdot,\cdot)\|_{(2p)^{'},(2q)^{'},\varrho_{x}^{-1}}
\end{align*}
which is finite for every $(s,x)\in Q_{\hat{T}}$ by Proposition
\ref{prop.p}. Now, let $u$ be a weak solution of the Cauchy
problem (\ref{eq4.6}) with $\bar f=0,\varphi=0$. Then
$E_{s,x}\int_s^T|\nabla u(\theta,X_\theta)|^2\,d\theta\le
C\varrho^{-2}(x)\|f\|^2_{2p,2q,\varrho}$ by \cite[Theorem
4.1]{Roz.Div}, and $|u(s,x)|\le
\varrho^{-1}(x)C\|f\|_{2p,2q,\varrho}$ by Proposition
\ref{prop.pp}. Moreover, from \cite{Roz.Div} we know that
\[
u(t,X_t)=\int_t^Ta^{-1}\bar{f}(\theta,X_\theta)\,d^*X_\theta
-\int_t^T\nabla u(\theta,X_\theta)\,dM_{s,\theta},\quad t\in[s,T],\quad
P_{s,x}\mbox{-a.s.}
\]
By the above estimates and  Doob's ${\mathbb{L}}_2$-inequality,
\begin{align*}
E_{s,x}\left|\int_s^Tf(t,X_t)\,d^*X_t\right|^2 &\le
C(E_{s,x}\sup_{s\le t\le T}\varrho^{-2}(X_{t})+\varrho^{-2}(x))
\|f\|^2_{2p,2q,\varrho}\\ &\le
2C\varrho^{-2}(x)\|f\|^2_{2p,2q,\varrho}\,,
\end{align*}
the last inequality being a consequence of Corollary
\ref{moments}.
\\
(iii) Follows from Proposition \ref{prop3.1}.
\\
(iv) The fact that integrals with respect to backward and forward
martingale are well defined follows directly from (iii). Modifying
slightly \cite[Lemma 5.2]{Roz.Dirichlet} to the case of
time-inhomogeneous diffusions we have for $\alpha>0$,
\begin{align}\label{nq7.a}
&\left| E_{s,x}\int_{s}^{T} |\bar{f}(t,X_t)|\,
d|\alpha^{s,x}_{s,\cdot}|_t)\right|^2 \le\left|\int_{Q_{sT}}|f|
|\frac{\partial p}{\partial y_j}|(s,x,t,y)\,dt\,dy\right|^{2}\nonumber\\
&\le CE_{s,x}\int_{s}^{T}(t-s)^{-\alpha} |\bar
f(t,X_t)|^2\,dt\times \sum^d_{i,j=1}\int_{Q_{sT}}(t-s)^{\alpha}
p^{-1}a_{ij}\frac{\partial p}{\partial y_i}
\frac{\partial p}{\partial y_j}(s,x,t,y)\,dt\,dy\nonumber\\
&\le CE_{s,x}\int_{s}^{T}(t-s)^{-\alpha}|\bar f(t,X_t)|^2\,dt.
\end{align}
Multiplying this inequality by $\varrho$ and using the fact that
the measure $\varrho\,dm$ is finite on ${\mathbb{R}}^{d}$ we
obtain by Jensen's inequality that
\begin{align*}
&I\equiv\int_{Q_T}\left(E_{s,x}\int_{s}^{T}
|\bar{f}(t,X_t)|\,d|\alpha^{s,x}_{s,\cdot}|_t\right)
\varrho(x)\,dx\,ds\\ &\quad \le
C\left(\int_{0}^{T}\!\!\int_{s}^{T}(t-s)^{-\alpha}
\left(\int_{{\mathbb{R}}^{d}}
E_{s,x}|f(t,X_t)|^2\varrho^{2}(x)\,dx\right)\,dt\,ds\right)^{1/2}.
\end{align*}
Write $r(t)=\|f(t)\|_{2,\varrho}^{2}$. From the above with
$\alpha=1/2$  and (iii) we get
\begin{align*}
I^2&\le C\int_{0}^{T}\!\!\int_{s}^{T}(t-s)^{-1/2}r(t)\,dt \,ds
=C\int_{0}^{T}\left(\int_{0}^t(t-s)^{-1/2}\,ds\right)r(t)\,dt\\
&\le CT^{1/2}\int_{0}^{T}r(t)\,dt =CT^{1/2}\|\bar
f\|^2_{2,\varrho,T}\,,
\end{align*}
which proves the result. The second assertion is a direct
consequence of  (iii) and (\ref{eq4.3}).
\end{proof}

\end{document}